\documentclass {article}
\usepackage{authblk}
\usepackage{etex}
\usepackage[utf8]{inputenc}
\usepackage[english]{babel}
\usepackage{amsmath}
\usepackage{amsthm}
\usepackage{amssymb}
\usepackage{sectsty}
\usepackage{titlesec}
\usepackage{color}
\usepackage[color,matrix,arrow]{xy}
\usepackage{amsgen}
\usepackage{amstext}
\usepackage{amsbsy}
\usepackage{amsopn}
\usepackage{amsfonts}
\usepackage{eepic}
\usepackage{graphicx}
\usepackage{epsf}
\usepackage{pstricks}
\usepackage{enumerate}
\usepackage{eqnarray}
\usepackage{faktor}
\xyoption{all}

\usepackage{pgf}
\usepackage{tikz}

\newcommand{\footrecall}[1]{%
}

\titleformat*{\section}{\large\bfseries}
\titleformat*{\subsection}{\normalsize \bfseries}

\newtheorem{theorem}{Theorem}[section]
\newtheorem{lemma}[theorem]{Lemma}

\newtheorem{remark}[theorem]{Remark}
\newtheorem{proposition}[theorem]{Proposition}
\newtheorem{corollary}[theorem]{Corollary}

\theoremstyle{definition}
\newtheorem{definition}[theorem]{Definition}

\newcommand{\N}{\mathbb{N}}
\newcommand{\Z}{\mathbb{Z}}

\newcommand{\Fix}{\text{Fix}}
\newcommand{\Sing}{\text{Sing}}
\newcommand{\Per}{\text{Per}}
\newcommand{\Ker}{\text{Ker}}
\newcommand{\End}{\text{End}}
\newcommand{\Mono}{\text{Mon}}
\newcommand{\rank}{\text{rk}}
\newcommand{\sol}{\text{Sol}}
\newcommand{\Aut}{\text{Aut}}
\newcommand{\Reg}{\text{Reg}}
\newcommand{\mc}{\mathcal}

\title{On endomorphisms of the direct product of two free groups}

\author{Andr\'e Carvalho \thanks{andrecruzcarvalho@gmail.com}}
\affil{Centre of Mathematics, University of Porto, R. Campo Alegre, 4169-007 Porto, Portugal}

\begin{document}

\maketitle

\begin{abstract}
  We describe the endomorphisms of the direct product of two free groups of finite rank and obtain conditions for which the subgroup of fixed points is finitely generated and we do the same for periodic points. We also describe the automorphisms of such a group and solve the three Whitehead problems for automorphisms, monomorphisms and endomorphisms for groups in this class. We also obtain conditions for an endomorphism to be uniformly continuous when a suitable metric is considered and study infinite fixed points of the extension to the completion of the group.\\
  

\end{abstract}

\section{Introduction}
The dynamical study of endomorphisms of groups started with the (independent) work of Gersten \cite{[Ger87]} and Cooper \cite{[Coo87]}, using respectively graph-theoretic and topological approaches. They proved that the subgroup of fixed points $\Fix(\varphi)$ of some fixed automorphism $\varphi$ of $F_n$ is always finitely generated, and Cooper succeeded on classifying from the dynamical viewpoint the fixed points of the continuous extension of $\varphi$ to the boundary of $F_n$. Bestvina and Handel subsequently developed the theory of train tracks to prove that $\Fix(\varphi)$ has rank at most $n$ in \cite{[BH92]}. The problem of computing a basis for $\Fix(\varphi)$ had a tribulated history and was finally settled by Bogopolski and Maslakova in 2016 in \cite{[BM16]}.

This line of research extended early to wider classes of groups. For instance, Paulin proved in 1989 that the subgroup of fixed points of an automorphism of a hyperbolic group is finitely generated \cite{[Pau89]}. Fixed points were also studied for right-angled Artin groups \cite{[RSS13]} and lamplighter groups \cite{[MS18]}.

In this paper we will study endomorphisms of the direct product of two free groups $F_n\times F_m$, with $m,n\geq 2$. We remark that the endomorphisms of free-abelian times free groups $\Z^m\times F_n$ are already studied in \cite{[DV13]} and their extension to the completion is studied in \cite{[Car20]}. Although the class of free groups is well known, some problems in the product $F_n\times F_m$ are not easily reduced to problems in each factor. In particular, when endomorphisms (and automorphisms) are considered, we have that many endomorphisms are not obtained by applying an endomorphism of $F_n$ to the first component and one of $F_m$ to the second, so some problems arise when the dynamics of an endomorphism in the product is considered. In fact, this class, which is a subclass of the right-angled Artin groups, has some wild subgroups. For instance, there is a finitely generated subgroup $H$ of $F_2\times F_2$ for which it is undecidable if a given element $g\in F_2\times F_2$ belongs to $H$ (see \cite{[Mih58]}). Surprisingly, endomorphisms of $F_n\times F_m$ have a simple description.
 
The paper is organized as follows. In Section \ref{secpreliminares}, we present some preliminaries on groups and endomorphisms. In Section \ref{secendautos}, we describe the endomorphisms (and automorphisms) of $F_n\times F_m$ and prove that the Whitehead problems for endomorphisms, monomorphisms and automorphisms are solvable. We find conditions for an endomorphism to have finitely generated fixed subgroup in Section \ref{secfixed}. In Section \ref{secperiodic}, we do the same for the subgroup of periodic points. Finally, in Section \ref{seccompletamento}, we describe the uniformly continuous endomorphisms when a suitable metric is considered and prove some dynamical properties of the extension of an endomorphism to the completion.

\section{Preliminaries}
\label{secpreliminares}
Given a group $G$, we denote by $\End(G)$ (resp. $\Mono(G),\Aut(G)$) the endomorphism monoid (resp. monomorphism monoid, automorphism group) of $G$. Given an endomorphism $\varphi$, let $\Fix(\varphi)=\{x\in G\mid x\varphi=x \}$  and $$\Per(\varphi)=\bigcup\limits_{k\geq 1} \Fix(\varphi^n)$$ be the \emph{fixed subgroup of $\varphi$} and the \emph{periodic subgroup of $\varphi$}, respectively.

We denote by $F_n, F_m$ the free groups of rank $n$ and $m$ and take  $A=\{a_1,\ldots,a_n\}$ and $B=\{b_1,\ldots,b_m\}$ to be the alphabets of $F_n$ and $F_m$, respectively. We will always assume that $m,n>1.$

We want to provide conditions for the fixed subgroup to be finitely generated. For that, we will use the following theorem by Anisimov and Seifert. Given a finite alphabet $A$ and a homomorphism $\pi:A^*\to G$, the \emph{rational subsets} of $G$ are the images through $\pi$ of rational $A$-languages. The set of all such subsets is denoted by $\text{Rat } G$
\begin{theorem}[\cite{[Ber79]}, Theorem III.2.7]
\label{AnisimovSeifert}
Let $H$ be a subgroup of a group $G$. Then $H\in \text{Rat } G$ if and only if $H$ is finitely generated.
\end{theorem}

Given two reduced words $u$ and $v$ in the free group, we write $u\wedge v$ to denote the longest common prefix of $u$ and $v$.

The prefix metric on a free group is defined by
$$d(u,v)=\begin{cases}
2^{-|u\wedge v|} \text{ if $u\neq v$}\\
0 \text{ otherwise}
\end{cases}.$$

The prefix metric on a free group is in fact an ultrametric and its completion $(\hat{F}_n,\hat{d})$ is a compact space which can be described as the set of all finite and infinite reduced words on the alphabet $A \cup A^{-1}$. We will denote by $\partial F_n$ the set consisting of only the infinite words and call it the \emph{boundary} of $F_n$.

We consider the direct product of two free groups $F_n\times F_m$ endowed with the product metric given by taking the prefix metric in each component, i.e.,
 $$d((x,y),(z,w))=\max\{d(x,z), d(y,w)\}.$$

This metric is also an ultrametric and $\widehat { F_n\times F_m}$ is homeomorphic to $\hat { F_n}\times \hat {F_m}$ by uniqueness of the completion (Theorem 24.4 in \cite{[SW70]}).

Also, throughout this paper we will often refer to a matrix and the linear application it defines indistinctively when no confusion arises.

\section{Endomorphisms and Automorphisms}
 \label{secendautos}
 
 In this section, we describe the endomorphisms and automorphisms of $F_n\times F_m$ and solve the Whitehead problems for $F_n\times F_m$. As usual, given $n\in \N$, we will denote the set $\{1,\ldots, n\}$ by $[n]$.
 \subsection{Endomorphisms}
Consider an endomorphism $\varphi: F_n\times F_m\to F_n\times F_m$ defined by $(a_i,1)\mapsto (x_i,y_i)$ and $(1,b_j)\mapsto (z_j,w_j)$ for $i\in [n]$ and $j\in [m].$ We define $X=\{x_i\mid i\in[n]\}$, $Y=\{y_i\mid i\in[n]\}$, $Z=\{z_j\mid j\in[m]\}$ and $W=\{w_j\mid j\in[m]\}$. We say that these sets are \emph{trivial} if they are singletons containing only the empty word and \emph{nontrivial} otherwise. 

For $\varphi$ to be well defined, we must have that $x_iz_j=z_jx_i$ and $y_iw_j=w_jy_i$, for every $i\in [n]$ and $j\in [m]$. By Proposition 1.3.2 in \cite{[Lot83]}, we have that two words in a free group commute if and only if they are powers of the same word, thus, for every $(i,j)\in [n]\times [m]$, we have
\begin{align}
\label{potpalavraxz}
x_i=1 \, \vee \, z_j=1 \, \vee\, \left(\exists\, u\in F_n\setminus\{1\} \,\exists\,p_i,r_j\in \mathbb Z\setminus\{0\}: x_i=u^{p_i }\wedge z_j=u^{r_j}\right)
\end{align}
and a similar condition holds for $y_i$ and $w_j$: 
\begin{align}
\label{potpalavrayw}
y_i=1 \, \vee \, w_j=1 \, \vee\, \left(\exists\, v\in F_m\setminus\{1\} \,\exists\, q_i,s_j\in \mathbb Z\setminus\{0\}: y_i=v^{q_i }\wedge w_j=v^{s_j}\right)
\end{align}
\begin{lemma}
Suppose $X$ and $Z$  are nontrivial. Then, there is some $1\neq u\in F_n$ such that $X\cup Z\subseteq \{u^k\mid  k\in \Z\}$. Similarly, if $Y$ and $W$  are nontrivial, then, there is some $1\neq v\in F_n$ such that $Y\cup W\subseteq \{v^k\mid k\in \Z\}$.
\end{lemma}
\begin{proof} 
Just consider $\tilde X=X\setminus \{1\}$ and $\tilde Z=Z\setminus\{1\}$ which are both finite sets and apply (\ref{potpalavraxz}) or  (\ref{potpalavrayw}) . 
 \end{proof}

So, we will consider several different cases:
\begin{enumerate}[(I)]
\item All sets $X,Y,Z$ and $W$ are nontrivial
\item $X$ is the only trivial set
\item $Y$ is the only trivial set
\item $X$ and $Y$ are the only trivial sets
\item $X$ and $Z$ are the only trivial sets
\item $Y$ and $Z$ are trivial sets
\item $X$ and $W$ are trivial sets.
\end{enumerate}

These cases are sufficient since every other case is analogous to one of the above, by swapping order of the factors (notice we are not assuming any relation between $m$ and $n$). Indeed, if $W$ is the only trivial set, we reduce to the second case; if $Z$ is the only trivial set, we reduce to the third case;  if $Z$ and $W$ are the only trivial sets, we reduce to the fourth case; if $Y$ and $W$ are the only trivial sets, we reduce to the fifth case. If three or more sets are trivial, then we must fall into one of the two last cases.

We define $P=\{p_i\in\Z\mid i\in[n]\}$, $Q=\{q_i\in\Z\mid i\in[n]\}$, $R=\{r_j\in\Z\mid j\in[m]\}$ and $S=\{s_j\in\Z\mid j\in[m]\}$, where $p_i,q_i,r_j,s_j$ are the numbers from (\ref{potpalavraxz}) and (\ref{potpalavrayw}).
 So, matching the numeration above, $\varphi$ must have one of the following forms
\begin{enumerate}[(I)]
\item $(a_i,1)\mapsto (u^{p_i}, v^{q_i})$ and $(1,b_j)\mapsto (u^{r_j}, v^{s_j})$, for some $1\neq u \in F_n$, $1\neq v \in F_m$ and integers $p_i,q_i, r_j, s_j\in \mathbb Z$ for $(i,j)\in [n]\times [m]$, such that $P,Q,R,S\neq \{0\}$. 
\item $(a_i,1)\mapsto (1, v^{q_i})$ and $(1,b_j)\mapsto (z_j, v^{s_j})$, for some $1\neq v \in F_m$, $z_j\in F_n$ and integers $q_i,  s_j\in \mathbb Z$ for $(i,j)\in [n]\times [m]$,  such that $Q,S\neq \{0\}$ and $Z\neq \{1\}$. We will denote by $\phi$ the homomorphism from $F_m$ to $F_n$ mapping $b_j$ to $z_j$, $j\in[m]$.
\item $(a_i,1)\mapsto (u^{p_i}, 1)$ and $(1,b_j)\mapsto (u^{r_j}, w_j)$, for some $1\neq u \in F_n$, $w_j\in F_n$ and integers $p_i,  r_j\in \mathbb Z$ for $(i,j)\in [n]\times [m]$,  such that $P,R\neq \{0\}$ and $W\neq \{1\}$. We will denote by $\phi$ the endomorphism of $F_m$ mapping $b_j$ to $w_j$, $j\in[m]$.
\item$(a_i,1)\mapsto (1,1)$ and $(1,b_j)\mapsto (z_j,w_j)$, for some $(z_j,w_j) \in F_n\times F_m$ such that $Z, W\neq \{1\}$. We will denote the component mappings $\phi:F_m\to F_n$ and $\psi\in \End(F_m)$, defined by $b_j\mapsto z_j$ and $b_j\mapsto w_j$, $j\in [m]$, respectively.
\item $(a_i,1)\mapsto (1,v^{q_i})$ and $(1,b_j)\mapsto (1, v^{s_j})$, for some $1\neq v \in F_m$, and integers $q_i,  s_j\in \mathbb Z$ for $(i,j)\in [n]\times [m]$, such that $Q,S\neq \{0\}$.
\item $(a_i,1)\mapsto (x_i,1)$ and $(1,b_j)\mapsto (1,w_j)$, for some $x_i\in F_n$ and $w_j\in F_n$ for $(i,j)\in [n]\times [m]$.  We  will denote the component mappings $\phi\in \End(F_n)$ and $\psi\in \End(F_m)$, defined by $a_i\mapsto x_i$, $i\in [n]$ and $b_j\mapsto w_j$, $j\in[m]$, respectively.
\item $(a_i,1)\mapsto (1,y_i)$ and $(1,b_j)\mapsto (z_j,1)$, for some $y_i\in F_m$ and $z_j\in F_n$ for $(i,j)\in [n]\times [m]$. We  will denote the component mappings $\phi:F_n\to F_m$ and $\psi: F_m\to F_n$, defined by $a_i\mapsto y_i$, $i\in [n]$ and $b_j\mapsto z_j$, $j\in[m]$, respectively.\\
\end{enumerate}

 For $(i,j)\in[n]\times [m]$, we define   $\lambda_i:F_n\to \mathbb Z$ as the endomorphism given by $a_k\mapsto \delta_{ik}$ and $\tau_j:F_m\to \mathbb Z$ given by $b_k\mapsto \delta_{jk}$, where $\delta_{ij}$ is the Kronecker symbol.

 \subsection{Automorphisms}
We are now able to describe the automorphisms of $F_n\times F_m$.
\begin{proposition}
\label{monos}
An endomorphism $\varphi: F_n\times F_m$ is surjective if and only if it is of type VI or type VII such that the component mappings $\phi$ and $\psi$ are surjective. An endomorphism $\varphi: F_n\times F_m$ is injective if and only if it is of type VI or type VII such that the component mappings $\phi$ and $\psi$ are injective. 

\end{proposition}

\begin{proof} 
It is easy to check that type I, II and III endomorphisms are neither injective nor surjective.

A type IV endomorphism $\varphi$ is never surjective: we have that $F_n$ is in the kernel of $\varphi$ and there is no surjective homomorphism from $F_m$ to $F_n\times F_m$: $F_n\times F_m$ surjects to $\Z^{m+n}$ (through abelianization) which has rank $n+m$.
Also, it is never injective since for every $(x,y)\in F_n\times F_m$, we have that $(1,y)\varphi=(x,y)\varphi$.

Type V endomorphisms are clearly not surjective, since for every $1\neq u\in F_n$, we have that $(u,1)\not\in Im(\varphi)$.  Also, $(a_1a_2,1)\varphi=(1,v^{q_1+q_2})=(a_2a_1,1)\varphi$, so it is not injective.

It is obvious that a type VI endomorphism is surjective (resp. injective) if and only if both $\phi$ and $\psi$ are surjective (resp. injective) and the same holds for type VII endomorphisms.
\end{proof}

\begin{corollary}
\label{autos}
$\varphi\in \End(F_n\times F_m)$ is an automorphism if and only if one of the following holds:
\begin{enumerate}[i.]
\item $\varphi$ is a type VI endomorphism such that both $\phi$ and $\psi$ are automorphisms. 
\item $n=m$ and $\varphi$ is a type VII endomorphism such that both $\phi$ and $\psi$ are automorphisms. 
\end{enumerate}
\end{corollary}

\begin{corollary}
 Let $\theta_n\in \Aut(F_n\times F_n)$ be the involution defined by $(x,y)\mapsto (y,x).$ If $n\neq m$, then $$\Aut(F_n\times F_m)\cong \Aut(F_n)\times \Aut(F_m).$$ If $n=m$, then 
 $Aut(F_n\times F_n)$
is the semidirect product of $Aut(F_n)\times Aut(F_n)$ and $\langle \theta_n\rangle$.
\end{corollary}
\begin{proof} If $n\neq m$, the claim follows from  Corollary \ref{autos}. 

Suppose now that $n=m$ and denote by $\Aut_6(F_n\times F_n)$ the subgroup of $\Aut(F_n\times F_n)$ whose elements are type VI automorphisms and by $\Aut_7(F_n\times F_n)$ the subset of type VII automorphisms.
 It is easy to check that $\Aut_6(F_n\times F_n)\cong \Aut(F_n)\times \Aut(F_n)$ is a normal subgroup of $\Aut(F_n\times F_n)$ of index $2$, so, since $\theta_n$ has order $2$, we have that  $\Aut(F_n\times F_n)$ is the semidirect product of  $\Aut(F_n)\times \Aut(F_n)$ and $\langle \theta_n\rangle$.
 \end{proof}

\begin{corollary}
 $\Aut(F_n\times F_m)$ is finitely presented.\qed
\end{corollary}

\begin{corollary}
$F_n\times F_m$ is  hopfian but not cohopfian.
\end{corollary}

\begin{proof} If $\varphi\in \End(F_n\times F_m)$ is type VII and surjective, then $m=n$ and the component mappings are surjective endomorphisms of $F_n$. Since $F_n$ is hopfian, we have that the components are injective, so $\varphi$ is injective. If $\varphi$ is type VI, then the component mappings $\phi$ and $\psi$ are surjective endomorphisms of $F_n$ and $F_m$, respectively, thus injective and we have that $\varphi$ is injective.
 
It is not cohopfian since neither $F_n$ nor $F_m$ are cohopfian. Indeed take $\phi$ and $\psi$ injective and nonsurjective endomorphisms of $F_n$ and $F_m$ respectively and consider the type VI endomorphism $\varphi\in \End(F_n\times F_m)$ defined by taking $\phi$ and $\psi$ as the component mappings. We have that $\varphi$ is injective but not surjective. 
\end{proof}

\subsection{Whitehead Problems}
We start by stating the Whitehead problems for automorphisms, monomorphisms and endomorphisms.

\textbf{Whitehead problems ($WhP_a(G)$, $WhP_m(G)$, $WhP_e(G)$):} Given two elements $u,v\in G$, decide whether there exists an automorphism (resp. monomorphism, endomorphism) $\phi$ of $G$ such that $u\phi=v$; if so, find one by giving images of generators.

In \cite{[Day09]} and \cite{[Day14]}, the author generalizes the techniques from Whitehead and applies them to partially commutative groups, being able to solve $WhP_a(G)$, when $G$ is partially commutative, which obviously solves $WhP_a(F_n\times F_m)$. Also, since it is possible to solve equations in such groups (see \cite{[DJK16]}), $WhP_e(F_n\times F_m)$ can also be seen to be decidable.
The author is not aware of any previous result implying the solution of $WhP_m(F_n\times F_m)$. We will present an alternative way to solving the Whitehead problems using the classification of endomorphisms already achieved and the already known solution of the problems for free groups.

\begin{theorem} For $n\geq 2$
\label{Whfree}
\begin{enumerate}[(i)]
\item (\cite{[Whi36]}) $WhP_a(F_n)$ is solvable
\item (\cite{[CH10]}) $WhP_m(F_n)$ is solvable
\item (\cite{[Mak82]}) $WhP_e(F_n)$ is solvable
\end{enumerate}
\end{theorem}

\begin{corollary}
\label{Whdifranks}
Given $m,n\in \N$, free groups $F_n$, $F_m$ and elements $u\in F_n$, $v\in F_m$,  the problem of deciding whether there exists a homomorphism $\phi:F_n\to F_m$ such that $u\phi=v$ and in case it does, finding it, is decidable.
\end{corollary}
\begin{proof} If $n=m$, we have the Whitehead problem $WhP_e(F_n)$ for free groups, which is solvable. 

If $n<m$, then  consider $\psi:F_n\to F_m$ to be the natural embedding. Then, we can decide if we have an endomorphism $\varphi\in \End(F_m)$ such that $u\psi\varphi=v$. If it exists, then we take $\phi=\psi\varphi$.
If there is no such $\varphi$, then there is no $\phi:F_n\to F_m$ such that $u\phi=v$. Indeed, if there was
 some   $\phi:F_n\to F_m$ such that $u\phi=v$, defining  $\varphi\in End(F_m)$ to be the endomorphism that maps the first $n$ letters in the basis of $F_m$ through the endomorphism induced by $\phi$ and the $m-n$ extra elements in the basis of $F_m$ to $1$, we would have that  $\phi=\psi\varphi$.

 In case $n>m$, we proceed in a similar way, extending the codomain instead of the domain. So, take  $\psi:F_m\to F_n$ to be the natural embedding.  If there exists a homomorphism  $\phi:F_n\to F_m$ such that $u\phi=v$, then  $\phi\psi$ is an endomorphism of $F_n$ such that $u\phi\psi=v\psi$.
Now, suppose that there is an endomorphism  $\varphi\in End(F_n)$ such that $u\varphi=v\psi$. Consider $\theta:F_n\to F_m$ such that $\psi\theta=1$. It follows that, putting $\phi=\varphi\theta$, we have that $u\phi=v.$
\end{proof}

Following the proof for free groups monomorphisms in \cite{[CH10]} step by step, the same result follows when the free groups have different ranks.

\begin{corollary}
\label{monoranksdif}
Given $m,n\in \N$, free groups $F_n$, $F_m$ and elements $u\in F_n$, $v\in F_m$,  the problem of deciding whether there exists an injective homomorphism $\phi:F_n\to F_m$ mapping $u$ to $v$ is decidable.
\end{corollary}

We remark that, 
given $u\in F_n$, we can compute $U=\{k \mid \exists \alpha\in F_n : u=\alpha^k\}$.

\begin{proposition}
$WhP_a(F_n\times F_n)$, $WhP_m(F_n\times F_m)$
 and $WhP_e(F_n\times F_m)$ are solvable.
\end{proposition}
\begin{proof} Let $(x,y), (z,w)\in F_n\times F_m$. We want to check if there is an endomorphism (resp. monomorphism,
 automorphism) $\varphi$ such that $(x,y)\varphi=(z,w)$.

$WhP_a(F_n\times F_m)$ follows directly from Corollary \ref{autos} and Theorem \ref{Whfree} (i).

For $WhP_m(F_n\times F_m)$, we decide if there is a monomorphism of type VI and if not,  we check if there is a monomorphism of type VII. For type VI endomorphisms, it follows from Proposition \ref{monos} and   Theorem \ref{Whfree} (ii).  For type VII monomorphisms, we use Proposition \ref{monos} and Corollary \ref{monoranksdif}.

For $WhP_e(F_n\times F_m)$, we will check if there is a type I endomorphism such that $(x,y)\varphi=(z,w)$. If there is, we stop. If not, we check the existence of a type II endomorphism. If there is one, we stop. If not, we check the existence of endomorphisms of a type III endomorphism, and so on. 

So, we want to check if, given $(x,y),(z,w)\in F_n\times F_m$, there is a type I endomorphism $\varphi$ such that $(x,y)\varphi=(z,w)$, i.e., if there are $u\in F_n$, $v\in F_m$, $p_i\in \Z$, $q_i\in \Z$, $r_j\in \Z$ and $s_j\in \Z$, for $(i,j)\in[n]\times[m]$ such that 
$$\begin{cases}
z=u^{\sum\limits_{i\in[n]} \lambda_i(x)p_i+\sum\limits_{j\in[m]} \tau_j(y)r_j}\\
w=v^{\sum\limits_{i\in[n]} \lambda_i(x)q_i+\sum\limits_{j\in[m]} \tau_j(y)s_j}
\end{cases}$$

We compute the values of $\lambda_i(x)$, $\tau_j(y)$ for $(i,j)\in[n]\times[m]$, $Z=\{k\mid \exists \alpha\in F_n : z=\alpha^k\}$ and $W=\{k\mid \exists \alpha\in F_m : w=\alpha^k\}$. Then for every $(k,\ell)\in Z\times W$ we see if the linear Diophantine system 
$$\begin{cases}
k=\sum\limits_{i\in[n]} \lambda_i(x)p_i+\sum\limits_{j\in[m]} \tau_j(y)r_j\\
\ell=\sum\limits_{i\in[n]} \lambda_i(x)q_i+\sum\limits_{j\in[m]} \tau_j(y)s_j
\end{cases}$$
has a solution on $p_i,q_i, r_j,s_j$. If it does, then there is an endomorphism given by the solution of the system together with $u,v$ such that $u^k=z$ and $v^\ell=w$. If not, there is no endomorphism and we check the existence of a type II endomorphism.

To check the existence of such a type II endomorphism, we use Corollary  \ref{Whdifranks} to check if there is an endomorphism $\phi: F_m\to F_n$ such that $z=y\phi$. If there is not, then we stop and see if there is a type III endomorphism. If there is, we compute the values of $\lambda_i(x)$, $\tau_j(y)$ for $(i,j)\in[n]\times[m]$  and $W=\{k \mid \exists \alpha\in F_m : w=\alpha^k\}$. Then for every $k\in W$, we see if the linear Diophantine equation $$k=\sum\limits_{i\in[n]} \lambda_i(x)q_i+\sum\limits_{j\in[m]} \tau_j(y)s_j$$ has a solution on $q_i,s_j$.  If it does, then there is an endomorphism defined by $$(\alpha,\beta)\mapsto\left(\beta\phi,v^{\sum\limits_{i\in[n]} \lambda_i(\alpha)q_i+\sum\limits_{j\in[m]} \tau_j(\beta)s_j}\right),$$ for $v\in F_m$  such that $w=v^k$. If not, there is no endomorphism.

The type III case is entirely analogous to the type II.

Type IV reduces to Corollary  \ref{Whdifranks}  and Theorem \ref{Whfree} (iii).

Type V follows from decidability of the word problem in $F_n$ and the same argument as above for the second component.

Type VI follows from Theorem \ref{Whfree} (iii) for $F_n$ and $F_m$.

Type VII follows from Corollary \ref{Whdifranks}.
\end{proof}

\begin{remark}
In the multiple Whitehead problem, we are given two $k$-tuples of elements of $F_n\times F_m$, $(g_1,\ldots, g_k)$ and $(h_1,\ldots, h_k)$, and we want to decide whether there exists an endomorphism (or monomorphism or automorphism) that maps the $g_i$ to $h_i$, for all $i\in[k]$. We remark that, proceeding as above and using the corresponding known results for free groups, the multiple Whitehead problems for endomorphisms, monomorphisms, and automorphisms are also decidable in $F_n\times F_m$.
\end{remark}

\section{Fixed subgroup of an endomorphism}
\label{secfixed}
In this section, we aim at giving conditions for the fixed subgroup of an endomorphism to be finitely generated. We will deal with each type of endomorphisms one by one. We now state the main result from this section.

Given a Diophantine equation $\sum_{i=1}^k a_ix_i=0$, we denote its solution set (as a subset of $\Z^k$) by $\sol(\sum_{i=1}^k a_ix_i)$.  Also, we denote the abelianization map by $\rho:F_m\to \Z^m$.
\begin{theorem}
Let  $\varphi\in \End(F_n\times F_m)$. Then, $\Fix(\varphi)$ is finitely generated if and only if one of the following holds:
\begin{enumerate}
\item $\varphi$ is a Type I, II, IV, V, VI or VII endomorphism
\item $\varphi$ is a Type III endomorphism such that $\sum\limits_{i\in[n]} \lambda_i( u)p_i\neq1$
\item  $\varphi$ is a Type III endomorphism  such that $\sum\limits_{i\in[n]} \lambda_i( u)p_i=1$ and $$\rank\left((\Fix(\phi))\rho\cap\sol\left(\sum\limits_{j\in[m]} r_jx_j=0\right)\right)= \rank((\Fix(\phi))\rho),$$
\item $\varphi$ is a Type III endomorphism  such that $\sum\limits_{i\in[n]} \lambda_i( u)p_i=1$ and $\Fix(\phi)=1$
\item $\varphi$ is a Type III endomorphism  such that $\sum\limits_{i\in[n]} \lambda_i( u)p_i=1$, $\Fix(\phi)$ is cyclic, $(\Fix(\varphi))\rho\neq \mathbf0$ and $(\Fix(\phi))\rho\cap \sol\left(\sum\limits_{j\in[m]} r_jx_j=0\right)= \mathbf 0$
\end{enumerate}

\end{theorem}

\subsection{Type I endomorphisms}
We will study the fixed point subgroup of such an endomorphism. To start, we obviously have that a fixed point $(u,v)\in F_n\times F_m$ must be of the form $(u^a,v^b)\in F_n\times F_m$ for some $a,b\in\Z.$ 
For $i\in [n]$, we have that $\varphi$ is defined by 
$$(x,y)\mapsto\left(u^{\sum\limits_{i\in[n]} \lambda_i(x)p_i+\sum\limits_{j\in[m]} \tau_j(y)r_j},v^{\sum\limits_{i\in[n]} \lambda_i(x)q_i+\sum\limits_{j\in[m]} \tau_j(y)s_j}\right).$$

For $x\in F_n, y\in F_m$, we denote $\sum\limits_{i\in[n]} \lambda_i(x)p_i$ by $x^{P}$; $\sum\limits_{j\in[m]} \tau_j(y)r_j$ by $y^{R}$; $\sum\limits_{i\in[n]} \lambda_i(x)q_i$ by $x^{Q}$ and $\sum\limits_{j\in[n]} \tau_j(y)s_i$ by $y^{S}$. We will keep this notation throughout the paper.

So 
\begin{align*}
(u^a,v^b)\in \Fix(\varphi)&\Leftrightarrow  
\begin{cases}
a=\sum\limits_{i\in[n]} \lambda_i (u^a)p_i+\sum\limits_{j\in[m]} \tau_j( v^b)r_j\\
b=\sum\limits_{i\in[n]} \lambda_i( u^a )q_i+\sum\limits_{j\in[m]} \tau_j( v^b )s_j
\end{cases}\\
&\Leftrightarrow  
\begin{cases}
a=au^P+bv^R\\
b=au^Q+bv^S
\end{cases}\\
&\Leftrightarrow  
\begin{cases}
a\left(-1+u^P\right)+bv^R=0\\
au^Q+b\left(-1+v^S\right)=0
\end{cases}
\end{align*}
So, consider the matrix given by $$M_\varphi=\begin{bmatrix} -1+u^P & v^R\\
u^Q& -1+v^S
 \end{bmatrix}$$
and we have that $\Fix(\varphi)=\{(u^a,v^b)\in F_n\times F_m \mid (a,b)\in \Ker(M_\varphi)\}$. It is in fact isomorphic to $Ker(M_\varphi)$ as a subgroup of the free-abelian group $\Z^2$, which is finitely generated and has a computable basis.
In particular, if 
$\text{det}(M_\varphi)
\neq 0$,
we have that the fixed subgroup of $\varphi$ is trivial.

\subsection{Type II endomorphisms}
Consider a homomorphism $\phi:F_m\to F_n$ given by $b_j\mapsto z_j$, $j\in[m]$. Then, we have that $\varphi$ is defined by 
$$(x,y)\mapsto\left(y\phi,v^{\sum\limits_{i\in[n]} \lambda_i(x)q_i+\sum\limits_{j\in[m]} \tau_j(y)s_j}\right).$$
So a fixed point must be of the form $(v^b\phi,v^b)$, for some $b\in \Z$. 
We have that 
\begin{align*}
(v^b\phi,v^b)\in \Fix(\varphi)
&\Leftrightarrow  
b=\sum\limits_{i\in[n]} \lambda_i( v^b\phi)q_i+\sum\limits_{j\in[m]} \tau_j(  v^b)s_j\\
&\Leftrightarrow b\left(-1+ ( v\phi)^Q +v^S\right)=0
\end{align*}

So, if $ ( v\phi)^Q +v^S\neq 1$, we have that $\Fix(\varphi)$ is trivial and otherwise, $\Fix(\varphi)$ is given by $\{(v^b\phi,v^b)\in F_n\times F_m \mid b\in \mathbb Z\}\cong \Z.$

\subsection{Type III endomorphisms}
Consider $\phi\in \End(F_m)$ given by $b_j\mapsto w_j$, $j\in[m]$.
Then, we have that $\varphi$ is defined by 
$$(x,y)\mapsto\left(u^{\sum\limits_{i\in[n]} \lambda_i(x)p_i+\sum\limits_{j\in[m]} \tau_j(y)r_j},y\phi\right).$$
So a fixed point must be of the form $(u^a,y)\in F_n\times F_m$, for some $a\in \Z$. 
We have that 
\begin{align*}
(u^a,y)\in \Fix(\varphi)
&\Leftrightarrow  
a=\sum\limits_{i\in[n]} \lambda_i(u^a)p_i+\sum\limits_{j\in[m]} \tau_j(y)r_j \quad \wedge \quad y\in \Fix(\phi)\\
&\Leftrightarrow a\left(-1+ u^P\right) +y^R=0 \quad \wedge \quad y\in \Fix(\phi)
\end{align*}

So, if  $\sum\limits_{i\in[n]} \lambda_i( u)p_i\neq1$, then putting $$G=\left\{y\in \Fix(\phi)\;\bigg\lvert\; \left(1-u^P\right) \text{ divides } {y^R}\right\}\leq \Fix(\phi)$$
we get
$$\Fix(\varphi)=\left\{\left(u^{\frac{y^R}{1- u^P}},y\right)\mid y\in G \right\}\cong G.$$

We now prove that $G$ is always finitely generated. Set $$H=\left\{y\in F_m\;\bigg\lvert\; \left(1-u^P\right) \text{ divides } {y^R}\right\}\leq F_m.$$

 We have that $G=\Fix(\phi)\cap H$ and $\Fix(\phi)$ is finitely generated. 
 
The subgroup $H$ has index $|1-u^P|$, and so it is finitely generated. Indeed, it follows from the fact that $\tau_j$ is a homomorphism for every $j\in [m]$, that  $$(xy)^R=\sum_{j=1}^m \tau_j(xy)r_j=\sum_{j=1}^m \tau_j(x)r_j+ \sum_{j=1}^m \tau_j(y)r_j=x^R+y^R$$
 holds for all  $x,y\in F_m$
 and so $xH=\{z\in F_m\mid x\equiv_{1-u^P} z\}$.
 By Howson's Theorem (\cite{[How54]}), we have that $G$ is finitely generated.

If, on the other hand, $\sum\limits_{i\in[n]} \lambda_i( u)p_i=1$, putting $H=\left\{y\in \Fix(\phi)\mid \sum\limits_{j\in[m]} \tau_j(y)r_j=0\right\}\leq \Fix(\phi)$, we have that 
$$\Fix(\varphi)=\{(u^a,y)\mid y\in H, a\in\Z \}\cong \mathbb Z\times H$$ and it is finitely generated if and only if $H$ is finitely generated.  
Notice this might not be the case. Indeed, if $m=2$, $r_1=1$, $r_2=-1$ and $\phi=Id$, we have that $$H=\left\{y\in \Fix(\phi)\mid \sum\limits_{j\in[m]} \tau_j(y)r_j=0\right\}=\{y\in F_2\mid\tau_1(y)=\tau_2(y)\}.$$
The language it defines is not rational, so, by Theorem \ref{AnisimovSeifert}, $H$ is not finitely generated.

Observe that, since $\tau_j(uvu^{-1})=\tau_j(v)$ for all $u,v\in F_m$, $j\in[m]$, then $H$ is a normal subgroup of $\Fix(\phi)$ (not of $F_n$) and so it is finitely generated if and only if it is trivial or has finite index in $\Fix(\phi)$. Consider the restriction $\rho':\Fix(\phi)\to (\Fix(\phi))\rho$ of the abelianization map $\rho$. Since $x\rho=(\tau_1(x),\ldots, \tau_m(x))$, we have that $$H\rho'= (\Fix(\phi))\rho'\cap \sol\left(\sum\limits_{j\in[m]} r_jx_j=0\right).$$
Also, $H=H\rho'\rho'^{-1}$: the fact that $H\subseteq H\rho'\rho'^{-1}$ is obvious and we have that $H\rho'\rho'^{-1}\subseteq H$ because $\Ker(\rho')\subseteq H.$
By \cite[Lemma 3.2 (ii)]{[DV13]}, we have that $H=H\rho'\rho'^{-1}$ has finite index in $\Fix(\phi)$ if and only if $H\rho'$ has finite index in $(\Fix(\phi))\rho'$, i.e. if and only if condition 3. holds:
$$\rank\left((\Fix(\phi))\rho\cap\sol\left(\sum\limits_{j\in[m]} r_jx_j=0\right)\right)= \rank((\Fix(\phi))\rho).$$
Now, we will describe the cases where $H$ is trivial. Since $\Ker(\rho')\subseteq H$, if $H$ is trivial then $\rho'$ is injective. If $\rank(\Fix(\phi))\geq 2$, then $\rho'$ cannot be injective since the commutator of two free generators is mapped to $1$. Thus, $\rho'$ is injective if and only if $\Fix(\phi)$ is trivial or $\Fix(\phi)$ is cyclic not abelianizing to zero. 

Also, if $\Fix(\phi)$ is trivial (condition 4.), then obviously $H$ is trivial and if $\Fix(\phi)$ is cyclic not abelianizing to zero, then $H$ is trivial if and only if $H\rho'$ is trivial (condition 5.).

Although the set of reduced words in $H$ might not be a rational language, we can prove that it is always a context-free language constructing a pushdown automata recognizing $H'=\left\{y\in F_m\mid \sum\limits_{j\in[m]} \tau_j(y)r_j=0\right\}$. Since the words in $\Fix(\phi)$ form a rational language, and $H=H'\cap \Fix(\phi)$, it follows that $H$ is context-free. To do so, consider $\Gamma=\{z,X,X^{-1}\}$ to be the stack alphabet and $z$ to be the starting symbol. When $b_i$ is read, if $r_i$ is positive, then if the top symbol of the stack has $X$ or $z$, we add $X^{r_i}$ to the stack. If the stack has at least $r_{i}$ $X^{-1}$'s on the top, we remove them and if there are $k<r_i$ symbols $X^{-1}$ at the top, we remove them and add $X^{r_i-k}$. If $r_i$ is negative, we do as above switching $X$ and $X^{-1}.$ Whenever $b_i^{-1}$ is read, we do the same switching $X$ and $X^{-1}$. If we read the starting symbol and have nothing to add, we reach a final state. This way, only a finite amount of memory is required, so $H'$ is context-free.

\subsection{Type IV endomorphisms}
Consider $\phi:F_m\to F_n$ defined by $b_j\mapsto z_j$, $j\in[m]$ and $\psi\in End(F_m)$ given by $b_j\mapsto w_j$, $j\in[m]$. Then $(x,y)\varphi=(y\phi,y\psi)$. In particular, $\Fix(\varphi)=\{(y\phi,y)\in F_n\times F_m \mid y\in\Fix(\psi)\}$ which is finitely generated. In fact, it is isomorphic to $\Fix(\psi)$ which is finitely generated  and from \cite{[Mut21]} has a computable basis. 

\subsection{Type V endomorphisms}
In this case, we have that $\varphi$ is defined by 
$$(x,y)\mapsto\left(1,v^{\sum\limits_{i\in[n]} \lambda_i(x)q_i+\sum\limits_{j\in[m]} \tau_j(y)s_j}\right).$$
In particular, we get that a fixed point must be of the form $(1,v^b)$, for some $b\in \Z$. We have that 
$$(1,v^b)\in \Fix(\varphi)\Leftrightarrow b=b\sum\limits_{j\in [m]} \tau_j(v)s_j.$$
It follows that, if $v^S=1$, then $\Fix(\varphi)\cong \mathbb Z$, since it is given by $\{(1,v^b)\mid b\in\Z\}$. If $v^S\neq 1$, then $\Fix(\varphi)$ is trivial.

\subsection{Type VI endomorphisms}
This case reduces to the free group case, since we have that $\varphi$ is given by $(x,y)\mapsto (x\phi,y\psi)$, for some $\phi\in\End(F_n)$ and $\psi\in \End(F_m).$ It follows that $\Fix(\varphi)=\Fix(\phi)\times\Fix(\psi)$, and so it is finitely generated and   from \cite{[Mut21]} has a computable basis. 

\subsection{Type VII endomorphisms}
This case is also similar to the previous, since we have that $\varphi$ is given by $(x,y)\mapsto (y\psi,x\phi)$, for $\phi:F_n\to F_m$ defined by $a_i\mapsto y_i$, $i\in[n]$ and $\psi:F_m\to F_n$ defined by $b_j\mapsto z_j$, $j\in[m]$. So a fixed point must be $(x,y)$ such that $x=y\psi$ and $y=x\phi$ and so $x=x\phi\psi$ and $y=y\psi\phi$. So, $\Fix(\varphi)\subseteq \Fix(\phi\psi)\times \Fix(\psi\phi)$. Also, for $x\in \Fix(\phi\psi)$, we have that $(x,x\phi)\varphi=(x\phi\psi,x\phi)=(x,x\phi)$. Similarly, we have that $(y\psi,y)\varphi=(y\psi,y)$ for $y\in \Fix(\psi\phi)$. So, $\Fix(\varphi)=\{(x,x\phi)\mid x\in \Fix(\phi\psi)\}=\{(y\psi,y)\mid y\in \Fix(\psi\phi)\}\cong \Fix(\phi\psi)\cong \Fix(\psi\phi)$ and thus it is finitely generated and  from \cite{[Mut21]} has a computable basis. \\

We remark that this also yields a somewhat trivial proof that for homomorphisms $\phi:F_n\to F_m$ and $\psi:F_m\to F_n$, the restrictions of $\phi$ and $\psi$ are explicit (mutually inverse) isomorphisms between $\Fix(\phi\psi)$ and $\Fix(\psi\phi).$

Since all the calculations above are explicit, we obtain the following algorithmic corollary.
\begin{corollary}
There is an algorithm which decides whether the fixed subgroup of a given endomorphism $\varphi\in \End(F_n\times F_m)$ is finitely generated, and it computes a set of generators (recursively, in the infinite case).
\end{corollary}

In the case of free groups, all endomorphisms have a finitely generated fixed subgroup, where in the case of free-abelian times free groups not even all automorphisms have a finitely generated fixed subgroup (see \cite[Example 6.4]{[DV13]}). In this class of groups, we have a situation in between these two situations.
\begin{corollary}
Let  $\varphi\in \Aut(F_n\times F_m)$. Then, $\Fix(\varphi)$ is finitely generated.
\end{corollary}

Also, since the proof above is constructive, we can present a complete list of isomorphism classes of finitely generated fixed subgroups of endomorphisms of $F_n\times F_m$. We say that a subgroup $H$ of a free group is \emph{context-free} if the language of reduced words in $H$ is context-free.
\begin{corollary}
Let  $\varphi\in \End(F_n\times F_m)$. Then, one of the following happens:
\begin{enumerate}
\item $\Fix(\varphi)$ is trivial;
\item $\Fix(\varphi)$ is isomorphic to $\Z^2$;
\item $\Fix(\varphi)$ is isomorphic to $\Z\times F$, where $F$ is a free group;
\item $\Fix(\varphi)$ is a finitely generated free group;
\item $\Fix(\varphi)$ is a direct product of free groups;
\item $\Fix(\varphi)$ is context-free.
\end{enumerate}
\end{corollary}
\noindent\textit{Proof.} If $\varphi$ is a type I endomorphism, then the fixed subgroup must be a subgroup of $\Z^2$ and so it must be trivial, infinite cyclic, or isomorphic to $\Z^2$. For type II endomorphisms, $\Fix(\varphi)$ must be trivial or isomorphic to $\Z$. For type III endomorphisms such that $\sum\limits_{i\in[n]} \lambda_i( u)p_i\neq1$, $\Fix(\varphi)$ is a finitely generated (in fact, finite index) subgroup of $F_m$ and so it is a finitely generated free group. If $\sum\limits_{i\in[n]} \lambda_i( u)p_i=1$, then it is isomorphic to $\Z\times H$, where $H$ is either trivial, a subgroup of a free group (hence free) or context-free. For type IV endomorphisms, $\Fix(\varphi)$ is isomorphic to the fixed subgroup of an endomorphism of $F_m$. For type V endomorphisms, $\Fix(\varphi)$ must be either isomorphic to $\Z$ or trivial. For type VI endomorphisms, it must be a direct product of fixed subgroups of endomorphisms of free groups and finally, for type VII, it is isomorphic to a fixed subgroup of an endomorphism of $F_m$.
\qed

\section{Periodic points}
\label{secperiodic} 
The purpose of this section is to describe the cases where $\Per(\varphi)$ is finitely generated. We will proceed as in the previous section, dealing with each type of endomorphisms one by one. 
We start by recalling that, in the case of free groups, periodic points have their period bounded by a computable constant depending only on the rank of the free group (\cite{[MS02]},\cite{[Car22]}). We will denote by $P_n$ the constant that bounds the periods for endomorphisms of $F_n$.

We now state the main result from this section.

\begin{theorem}
\label{periodico}
Let  $\varphi\in \End(F_n\times F_m)$. Then, $\Per(\varphi)$ is finitely generated if and only if one of the following holds:
\begin{enumerate}
\item $\varphi$ is a Type I, II, IV, V, VI or VII endomorphism
\item $\varphi$ is a Type III endomorphism such that $\left|\sum\limits_{i\in[n]} \lambda_i( u)p_i\right|\neq1$
\item  $\varphi$ is a Type III endomorphism  such that $\left|\sum\limits_{i\in[n]} \lambda_i( u)p_i\right|=1$ and $\Fix(\varphi^{(2P_m)!})$ is finitely generated.
\end{enumerate}
\end{theorem}

\subsection{Type I endomorphisms}
Consider a type I endomorphism given by 
$$(x,y)\mapsto\left(u^{\sum\limits_{i\in[n]} \lambda_i(x)p_i+\sum\limits_{j\in[m]} \tau_j(y)r_j},v^{\sum\limits_{i\in[n]} \lambda_i(x)q_i+\sum\limits_{j\in[m]} \tau_j(y)s_j}\right).$$

We start to determine the orbit of a point $(x,y)\in F_n\times F_m$. Define sequences in $\Z$ by $a_1(x,y)=x^P+y^R$, $b_1(x,y)=x^Q+y^S$ and $a_{n+1}(x,y)=a_n(x,y)u^P+b_n(x,y)v^R$ and $b_{n+1}(x,y)=a_n(x,y) u^Q+b_n(x,y) v^S$. We have that, for every $k\in \N$, $(x,y)\varphi^k=(u^{a_k(x,y)},v^{b_k(x,y)})$. Indeed, for $k=1$, we have that $(x,y)\varphi=(u^{a_1(x,y)},v^{b_1(x,y)})$. Assume that $(x,y)\varphi^r=(u^{a_r(x,y)},v^{b_r(x,y)}),$ for every $r\leq k$. We have that 
\begin{align*}
(x,y)\varphi^{k+1}&=(x,y)\varphi^{k}\varphi=(u^{a_k(x,y)},v^{b_k(x,y)})\varphi\\&=(u^{a_k(x,y)u^P+b_k(x,y)v^R},v^{a_k(x,y)u^Q+b_k(x,y)v^S})\\&=(u^{a_{k+1}(x,y)},v^{b_{k+1}(x,y)}).
\end{align*}

Clearly, a periodic point must be of the form $(x,y)=(u^a, v^b)$. In this case, $x^P=au^P$, $y^R=bv^R$, $x^Q=au^Q$ and $y^S=bv^S$, so putting $a_0=a$, $b_0=b$, $a_{n+1}=a_nu^P+b_nv^R$ and $b_{n+1}=a_n u^Q+b_n v^S$, we have that $(u^a, v^b)\varphi^k=(u^{a_k},v^{b_k})$. Defining the matrix $$M_\varphi=\begin{bmatrix} u^P &  v^R \\u^Q&v^S
\end{bmatrix}$$
and denoting also by $M_\varphi$ the endomorphism of $\Z^2$ defined by the matrix, we have that $$\Per(\varphi)=\{(u^a,v^b)\in F_n\times F_m\mid (a,b)\in \Per(M_\varphi)\}\cong \Per(M_\varphi),$$
which is finitely generated, as every subgroup of $\Z^2$.
\subsection{Type II endomorphisms}
Consider a type II endomorphism defined by 
$$(x,y)\mapsto\left(y\phi,v^{\sum\limits_{i\in[n]} \lambda_i(x)q_i+\sum\limits_{j\in[m]} \tau_j(y)s_j}\right).$$
We want to study the orbit of a point $(x,y)\in F_n\times F_m$. We define a sequence of integers $(a_n)_n$ by $a_1(x,y)=x^Q+ y^S$, $a_2(x,y)=a_1(x,y) v^S +(y\phi)^Q$ and for $n>2$, $a_n(x,y)=a_{n-1}(x,y)v^S +a_{n-2}(x,y)(v\phi)^Q.$ We will now prove that, for $k\geq 2$, we have that  $(x,y)\varphi^k=(v^{a_{k-1}(x,y)}\phi,v^{a_k(x,y)})$. For $k=2$, we have that $(x,y)\varphi^2=(y\phi,v^{x^Q+ y^S})\varphi=(v^{a_1(x,y)}\phi, v^{({y\phi})^Q+a_1(x,y)v^S})=(v^{a_1(x,y)}\phi, v^{a_2(x,y)})$. Now, assume the statement holds for every $r\leq k$. We have that 
\begin{align*}
(x,y)\varphi^{k+1}&=(v^{a_{k-1}(x,y)}\phi,v^{a_k(x,y)})\varphi=(v^{a_k(x,y)}\phi,v^{a_{k-1}(x,y)(v\phi)^Q+a_k(x,y)v^S})\\ &=(v^{a_k(x,y)}\phi,v^{a_{k+1}(x,y)}).
\end{align*}

So, a periodic point  has the form $(v^b\phi,v^a)$, for $a,b\in \Z.$ 

Suppose that $v\phi\neq 1$ and consider $$H=\{(b,a)\in \Z^2\mid (v^b\phi,v^a)\in \Per(\varphi)\}.$$

Clearly, $H$ is a subgroup of $\Z^2$, thus finitely generated and $\Per(\varphi)\cong H.$

Now, if $v\phi=1$, then a periodic point must be of the form $(x,y)=(1,v^a)$ which is isomorphic to a subgroup of $\Z$, hence finitely generated.

\subsection{Type III endomorphisms}
As in the fixed point case, the subgroup of periodic points of a Type III endomorphism is more complex than the subgroup of periodic points of endomorphisms of other types. In particular, it might not be finitely generated, as seen in the fixed point case.

Clearly, a periodic point must be of the form $(x,y)=(u^a,w)\in F_n\times F_m$. We will prove by induction on $k$ that 
\begin{align}
\label{orbitt3}
(u^a,w)\varphi^k=\left(u^{a(u^P)^{k}+\sum\limits_{t=0}^{k-1} (w\phi^t)^R(u^P)^{k-t-1}},w\phi^k\right).
\end{align} 
For $k=1$, we have that $(u^a,w)\varphi=(u^{au^P+w^R},w\phi)$. Now, assume the statement holds for every $r\leq k$. We have that 
\begin{align*}
(u^a,w)\varphi^{k+1}=&(u^a,w)\varphi^k\varphi=\left(u^{a(u^P)^{k}+\sum\limits_{t=0}^{k-1} (w\phi^t)^R(u^P)^{k-t-1}},w\phi^k\right)\varphi\\
=&\left(u^{\left(a(u^P)^{k}+\sum\limits_{t=0}^{k-1} (w\phi^t)^R(u^P)^{k-t-1}\right)u^P+(w\phi^k)^R},w\phi^{k+1}\right)\\
=&\left(u^{a(u^P)^{k+1}+\sum\limits_{t=0}^{k} (w\phi^t)^R(u^P)^{k-t}},w\phi^{k+1}\right).
\end{align*}

If $u^P=1$, then put $$H=\left\{y\in \Per(\phi)\;\lvert\; \exists s\in \N :  \sum\limits_{t=0}^{s\pi_y-1} (y\phi^t)^R=0 \right\}.$$
Observe that for $s\in \N$, $y\in \Per(\phi)$, we have that $\sum\limits_{t=0}^{s\pi_y-1} (y\phi^t)^R=s\sum\limits_{t=0}^{\pi_y-1} (y\phi^t)^R$, so 
$$H=\left\{y\in \Per(\phi)\;\lvert\; \sum\limits_{t=0}^{\pi_y-1} (y\phi^t)^R=0 \right\},$$
where $\pi_y$ denotes the period of $y$.
We have that $\Per(\varphi)=\{(u^a,w) \mid w\in H, a\in \Z\}\cong \Z\times H$, and we know that $H$ might not be finitely generated, as seen in the fixed points case.

Moreover, it follows that if $(u^a,y)$ is periodic, then its period is equal to the period of $y$, which is bounded by $P_m$, and so, in this case, $\Per(\varphi)=\Fix(\varphi^{P_m!})=\Fix(\varphi^{(2P_m)!}).$

If $u^P=-1$, we will consider two cases. Let $y\in \Per(\phi)$ with even period $\pi_y$ and take $a\in \Z$ and $s\in \N$.
We have that $$(u^a,y)\varphi^{s\pi_y}=\left(u^{a+\sum\limits_{t=0}^{s\pi_y-1}(y\phi^t)^R(-1)^{t+1}},y\right).$$ So, $(u^a,y)\in \Per(\varphi)$ if and only if there is $s\in \N$, such that $\sum\limits_{t=0}^{s\pi_y-1}(y\phi^t)^R(-1)^{t+1}=0$. Since, $\sum\limits_{t=0}^{s\pi_y-1}(y\phi^t)^R(-1)^{t+1}=s\sum\limits_{t=0}^{\pi_y-1}(y\phi^t)^R(-1)^{t+1}$, then putting $$H=\left\{y\in \Per(\phi)\mid \pi_y\equiv_2 0\wedge \sum\limits_{t=0}^{\pi_y-1}(-1)^{t+1}(y\phi^t)^R=0\right\},$$
we have that $$(u^a,y)\in \Per(\varphi)\Leftrightarrow y\in H.$$

It follows that, in this case, the period of $(u^a,y)$ is the same as the period of $y$, which is bounded  by $P_m$, and so, in this case, $\Per(\varphi)=\Fix(\varphi^{P_m!})=\Fix(\varphi^{(2P_m)!}).$

Now, let $y\in \Per(\phi)$ with odd period $\pi_y$ and take $a\in \Z$ and $s\in \N$. Then $$(u^a,y)\varphi^{2\pi_y}=\left(u^{a+\sum\limits_{t=0}^{2\pi_y-1}(-1)^{t+1}(y\phi^t)^R},y\right)=(u^a,y),$$
and the same equality is not true with any other positive exponent smaller than $2\pi_y.$

So, $\Per(\varphi)=\{(u^a,y)\in F_n\times \Per(\phi)\mid y\in H\vee \pi_y\equiv_2 1\}.$ Notice that, if every periodic point of $\phi$ has odd period, then $\Per(\varphi)$ is finitely generated, since in that case, $\Per(\varphi)\cong \Z\times \Per(\phi).$ However, it might be the case where $\Per(\varphi)$ is not finitely generated. For example, let $m=2$ and $\phi\in \End(F_2)$ defined by $a\mapsto b$ and $b\mapsto a$. Then $\Per(\phi)=F_m$ and every point has even period. If $r_1=1$ and $r_2=-1$, then 
\begin{align*}
H&=\{y\in F_2\mid -y^R+(y\phi)^R=0\}\\&=\{y\in F_2\mid -\lambda_1(y)+\lambda_2(y)+\lambda_2(y)-\lambda_1(y)=0\}\\&=\{y\in F_2\mid \lambda_1(y)-\lambda_2(y)=0\}
\end{align*}
which is not finitely generated.

In this case, if $(u^a,y)$ is periodic, then its period is at most $2\pi_y$, which is bounded by $2P_m$, and so, in this case, $\Per(\varphi)=\Fix(\varphi^{(2P_m)!}).$

In case $|u^P|\neq 1$, we will show that $\Per(\varphi)$ is always finitely generated. To do so, we start by proving that if $\varphi$ is a type III endomorphism such that $|u^P|\neq 1$, then for every $k\in \N$, we have that $\varphi^k$  has finitely generated fixed subgroup. Indeed, 
\begin{align*}
(a_i,1)\varphi^k=\left(u^{p_i(u^P)^{k-1}},1\right) \quad \text{ and } \quad (1,b_j)\varphi^k=\left(u^{\sum\limits_{t=0}^{k-1} (b_j\phi^t)^R(u^P)^{k-t-1}},b_j\phi^k\right),
\end{align*}
so if $\varphi^k$ is a type III endomorphism, then $$\left|\sum\limits_{i\in[n]}\lambda_i(u)p_i(u^P)^{k-1}\right|=\left|(u^P)^{k}\right|\neq 1.$$
So, we know that $\Fix(\varphi^k)$ is finitely generated for every $k\in \N$. We will now show that there is a bound for the periods of periodic points of $\varphi$. To do that, we will show that the period of $(u^a,w)$ is the same as the period of $w$ through $\phi$, $\pi_w$, which is bounded by $P_m$. This suffices to prove that $\Per(\varphi)$ is finitely generated since, in this case $$\Per(\varphi)=\bigcup\limits_{k=1}^{P_m!}\Fix(\varphi^k),$$
which is a finite union of finitely generated groups.

We have that $$(u^a,w)\varphi^k=(u^a,w)\Rightarrow k=s\pi_w,\text{ for some $s\in \N$}.$$

By (\ref{orbitt3}) putting $$S=\left\{w\in \Per(\phi)\mid \exists s_w\in N: (1-(u^P)^{s_w\pi_w}) \text{ divides } \sum\limits_{t=0}^{s_w\pi_w-1} (w\phi^t)^R(u^P)^{s_w\pi_w-t-1} \right\},$$
we have that $$\Per(\varphi)=\left\{\left(u^{\frac{ \sum\limits_{t=0}^{s_w\pi_w-1} (w\phi^t)^R(u^P)^{s_w\pi_w-t-1}}{1-(u^P)^{s_w\pi_w}}},w\right)\mid w\in S\right\}.$$

Now, for $w\in S$, reordering the summands, we have that 

\begin{scriptsize}
\begin{align*}
&\sum\limits_{t=0}^{s_w\pi_w-1} (w\phi^t)^R(u^P)^{s_w\pi_w-t-1} =\\
&= w^R(u^P)^{s_w\pi_w-1}       &+\;\;\,&(w\phi)^R(u^P)^{s_w\pi_w-2}        &+\cdots     &+ \;\;\,&(w\phi^{\pi_w-1})^R&(u^P)^{(s_w-1)\pi_w}\\
&+ w^R(u^P)^{(s_w-1)\pi_w-1} &+\;\;\,&(w\phi)^R(u^P)^{(s_w-1)\pi_w-2}   &+\cdots     &+\;\;\,& (w\phi^{\pi_w-1})^R&(u^P)^{(s_w-2)\pi_w}\\
&&&\vdots&\vdots\\
&+ w^R(u^P)^{\pi_w-1} &+\;\;\,&(w\phi)^R(u^P)^{\pi_w-2}   &+\cdots     &+ \;\;\,&(w\phi^{\pi_w-1})^R&\\
&=w^R\sum\limits_{t=1}^{s_w} (u^P)^{t\pi_w-1}&+\;\;\,& (w\phi)^R\sum\limits_{t=1}^{s_w} (u^P)^{t\pi_w-2}&+\cdots     &+\;\;\,&(w\phi^{\pi_w-1})^R&\sum\limits_{t=1}^{s_w} (u^P)^{(t-1)\pi_w}\\
&=w^R(u^P)^{\pi_w-1}\frac{1-(u^P)^{s_w\pi_w}}{1-(u^P)^{\pi_w}} &+\;\;\,&(w\phi)^R(u^P)^{\pi_w-2}\frac{1-(u^P)^{s_w\pi_w}}{1-(u^P)^{\pi_w}}   &+\cdots     &+ \;\;\,&(w\phi^{\pi_w-1})^R&\frac{1-(u^P)^{s_w\pi_w}}{1-(u^P)^{\pi_w}}.
\end{align*}
\end{scriptsize}

So, we have that $$\sum\limits_{t=0}^{s_w\pi_w-1} (w\phi^t)^R(u^P)^{s_w\pi_w-t-1} =\frac{1-(u^P)^{s_w\pi_w}}{1-(u^P)^{\pi_w}}\left(\sum\limits_{t=0}^{\pi_w-1} (w\phi^t)^R(u^P)^{\pi_w-t-1}\right),$$
and  $1-(u^P)^{s_w\pi_w}$  divides  $\sum\limits_{t=0}^{s_w\pi_w-1} (w\phi^t)^R(u^P)^{s_w\pi_w-t-1}$ if and only if $1-(u^P)^{\pi_w}$ divides $\sum\limits_{t=0}^{\pi_w-1}(w\phi^t)^R (u^P)^{\pi_w-t-1}$, which is a property of $w$ that does not depend on $s_w$.

So, we proved that given a periodic point $(u^a,w)$, we have that the existence of $s\in \N$ such that $(u^a,w)\varphi^{s\pi_w}=(u^a,w)$ depends only on $w$ and if there is one, then $(u^a,w)\varphi^{s\pi_w}=(u^a,w)$ for every $s\in \N$. Thus, the period of $(u^a,w)$ is  equal to the period of $w$ through $\phi$, which is bounded.

\subsection{Type IV endomorphisms}
It is easy to see by induction that  $(u,v)\varphi^r=(v\psi^{r-1}\phi,v\psi^r)$, for every $r>0$. Indeed, it is true for $r=1$ and if we have that $(u,v)\varphi^r=(v\psi^{r-1}\phi,v\psi^r)$, then $(u,v)\varphi^{r+1}=(v\psi^{r-1}\phi,v\psi^r)\varphi=(v\psi^{r}\phi,v\psi^{r+1}).$
So $\Per(\varphi)=\{(y\psi^{\pi_y-1}\phi,y)\mid y\in \Per(\psi)\}\cong \Per(\psi)$, where $\pi_y$ denotes the period of $y$ through $\psi$.
\subsection{Type V endomorphisms}
Consider a type V endomorphism $\varphi$  defined by 
$$(x,y)\mapsto\left(1,v^{\sum\limits_{i\in[n]} \lambda_i(x)q_i+\sum\limits_{j\in[m]} \tau_j(y)s_j}\right).$$
We prove by induction that $(x,y)\varphi^k=\left(1,v^{(x^Q+y^S)(v^S)^{k-1}}\right)$. For $k=1$, it is clear it holds. Now, assume the statement holds for every $r\leq k$. We have that 
\begin{align*}
(x,y)\varphi^{k+1}=&(x,y)\varphi^k\varphi=\left(1,v^{(x^Q+y^S)(v^S)^{k-1}}\right)\varphi=\left(1,v^{(x^Q+y^S)(v^S)^{k}}\right).
\end{align*}

So, a periodic point must have the form $(1,v^b)$ for some $b\in \Z$. In this case, $$(1,v^b)\varphi^k=(1,v^{bv^S(v^S)^{k-1}})=\left(1,v^{b(v^S)^{k}}\right).$$

So, if $v^S\not\in\{-1,1\}$, then, $\Per(\varphi)$ is trivial. If $v^S\in\{-1,1\}$, then $\Per(\varphi)\cong \Z.$
\subsection{Type VI endomorphisms}
We have that $(x,y)\varphi^k=(x\phi^{k},y\psi^k)$, so $$\Per(\varphi)=\Per(\phi)\times \Per(\psi).$$
\subsection{Type VII endomorphisms}
Consider an endomorphism $\varphi$ defined by $(x,y)\varphi=(y\psi, x\phi)$, for $\phi:F_n\to F_m$ and $\psi:F_m\to F_n$. We will prove by induction that for $k\in \mathbb N_0$, $(x,y)\varphi^{2k}= (x(\phi\psi)^k,y(\psi\phi)^k)$ and $(x,y)\varphi^{2k+1}= (y(\psi\phi)^k\psi,x(\phi\psi)^k\phi)$. We have that $(x,y)\varphi=(y\psi,x\varphi)$. Now suppose the claim holds for $(x,y)\varphi^r$ for every $r$ up to some $k$.
We have that \begin{align*}(x,y)\varphi^{k+1}&=(x,y)\varphi^k\varphi=
\begin{cases} (x(\phi\psi)^{\frac k2},y(\psi\phi)^{\frac k2})\varphi, \quad &\text{if $k$ is even} \\ (y(\psi\phi)^{\frac{k-1}2}\psi,x(\phi\psi)^{\frac{k-1}2}\phi)\varphi, \quad &\text{if $k$ is odd} 
\end{cases}\\
&=\begin{cases} (y(\psi\phi)^{\frac k2}\psi, x(\phi\psi)^{\frac k2}\phi), \quad &\text{if $k+1$ is odd} \\ (x(\phi\psi)^{\frac{k+1}2}, y(\psi\phi)^{\frac{k+1}2}), \quad &\text{if $k+1$ is even} 
\end{cases}
\end{align*}

So, the periodic points are given by $\Per(\phi\psi)\times\Per(\psi\phi)$, which is finitely generated.

\begin{corollary}
 There is an algorithm which decides whether the periodic subgroup of a given endomorphism of $F_n\times F_m$ is finitely generated.
\end{corollary}

Notice that if $\Per(\varphi)$ is finitely generated, then the periods must be bounded above by the least common multiple of the periods of the generators. The following corollary is immediate from the proof of Theorem \ref{periodico}.
\begin{corollary}
Let $\varphi\in \End(F_n\times F_m)$. There is a constant $C_\varphi \in \N$ such that $\Per(\varphi)=\Fix(\varphi^{C_\varphi})$.
\end{corollary}

\section{Dynamics of infinite points}
\label{seccompletamento}
For $i=1,2$, denote by $\pi_i$ the projections in the $i$-th component and by $\varphi_i$ the homomorphism $\varphi\pi_i$.
We have that $\varphi$ is uniformly continuous if and only if both $\varphi_1$ and $\varphi_2$ are. Also, recall the definition of the sets $X,Y,Z,W$ introduced in Section \ref{secfixed}. We will consider $F_n\times F_m$ endowed with the product metric given by taking the prefix metric in each component. Check Section \ref{secpreliminares} for definitions.

If $\varphi$ is uniformly continuous, then it admits a unique continuous extension to the completion $\hat\varphi:\widehat{F_n\times F_m}\to\widehat{F_n\times F_m}$.  The dynamics of this extension has been studied for several classes of groups including the free groups (see, for example, \cite{[Coo87]},\cite{[GJLL98]},\cite{[LL08]}), virtually free groups (see \cite{[Sil13]}) and free-abelian times free groups (see \cite{[Car20]}).

A fixed point $\alpha\in \partial(F_n\times F_m)$ is said to be \emph{singular} if $\alpha$ belongs to the topological closure $(\Fix(\varphi))^c$ of $\Fix(\varphi)$. Otherwise $\alpha$ is said to be \emph{regular}. We denote by $\Sing(\hat\varphi)$ (resp. $\Reg(\hat\varphi)$) the set of all singular (resp. regular) infinite fixed points of $\hat\varphi$.

\begin{definition}
\label{defatrep}
Let $\varphi\in \Aut(F_n\times F_m).$ An infinite fixed point $\alpha\in \Fix(\hat\varphi)$ is 
\begin{itemize}
\item an \emph{attractor} if $$\exists \varepsilon >0\,\,\forall\beta\in\widehat{F_n \times F_m} (d(\alpha,\beta)<\varepsilon\Rightarrow \lim\limits_{n\to +\infty} \beta\hat\varphi^n=\alpha)$$
\item a \emph{repeller} if   $$\exists \varepsilon >0\,\,\forall\beta\in\widehat{F_n \times F_m}(d(\alpha,\beta)<\varepsilon\Rightarrow \lim\limits_{n\to +\infty} \beta\hat\varphi^{-n}=\alpha)$$
\end{itemize}
\end{definition}

In this section, we obtain conditions for an endomorphism of $F_n\times F_m$ to be uniformly continuous and see how some questions typically studied can be reduced to the free group case.

\begin{proposition}
Let $\varphi\in End(F_n\times F_m)$. The following conditions are equivalent:
\begin{enumerate}
\item $\varphi$ is uniformly continuous.
\item $\varphi$ is of type IV, VI or VII with uniformly continuous components functions $\phi$ and $\psi$.
\end{enumerate}
\end{proposition}
\begin{proof}
It is obvious that $2.\Rightarrow 1.$

$1.\Rightarrow 2.$ Suppose that $\varphi$ is uniformly continuous. We start by proving that $1\in X\Rightarrow X=1$. Indeed, if $x_i=1$ for some $i\in [n]$ and $x_j\neq 1$ for some $j\in[n]$, then for every $\delta>0$, we have  
$$d\left(\left(a_i^{\lceil log_2(\frac 1\delta)\rceil},1\right),\left(a_i^{\lceil log_2(\frac 1\delta)\rceil}a_j,1\right)\right)<\delta$$ and $$d\left(\left(a_i^{\lceil log_2(\frac 1\delta)\rceil},1\right)\varphi_1,\left(a_i^{\lceil log_2(\frac 1\delta)\rceil}a_j,1\right)\varphi_1\right)=d(1,x_j)=1.$$ The same holds for $Y,Z$ and $W$.

We now prove that we can't have $Y\neq \{1\}$ and $ W\neq \{1\}$ (types I, II and V). Suppose we do. Then $1\not\in Y\cup W.$  We know that there is some word $1\neq v\in F_m$ and nonzero exponents $q_i$ and $s_j$ such that $y_i=v^{q_i}$ and $w_j=v^{s_j}$. Take $\delta>0$ and set $k=\lceil log_2(\frac 1\delta)\rceil$. We have that $$d\left((a_i^{ks_j},b_j^{-kq_i}),(a_i^{ks_j+1},b_j^{-kq_i})\right)<\delta$$ and $$d\left((a_i^{ks_j},b_j^{-kq_i})\varphi_2,(a_i^{ks_j+1},b_j^{-kq_i})\varphi_2\right)=d(1,v^{q_i})=1.$$ The same argument shows that we can't have $X\neq \{1\}$ and $ Z\neq \{1\}$, so type III is also done.

It is obvious that an endomorphism of type IV, VI or VII is uniformly continuous if and only if both $\phi$ and $\psi$ are.
 \end{proof}

Following the proof in \cite{[CS09b]} step by step, where the following lemma is proved for endomorphisms of a free group, we can actually prove the same result for homomorphisms of free groups with different ranks.
\begin{lemma}
Let $\phi:F_n\to F_m$ be a homomorphism. Then $\phi$ is uniformly continuous if and only if it is either trivial or injective.
\end{lemma}

Definition \ref{defatrep} is only presented for automorphisms because the definition of a repeller assumes the existence of an inverse. We will also consider the definition of an attractor in the cases where $\varphi$ is not an automorphism but still admits an extension to the completion.

\subsection{Type IV uniformly continuous endomorphisms}
Let $\varphi:F_n\times F_m\to F_n\times F_m$ be a type IV uniformly continuous endomorphism.  If both component mappings $\phi$ and $\psi$ are trivial, then $\varphi$ is trivial. If $\phi$ is trivial and $\psi$ is injective, then $(u,v)\varphi^r=(1,v\psi^r)$ and $\varphi$ has essentially the same dynamical behaviour as $\psi$. If $\psi$ is trivial and $\phi$ is injective, then $(u,v)\varphi=(v\phi,1)$ and $(u,v)\varphi^r=(1,1)$, for $r>1$. 

So, take injective $\phi\in F_m\to F_n$ and $\psi\in End(F_m)$ such that $ (u,v)\mapsto (v\phi,v\psi)$. We have that  $(u,v)\varphi^r=(v\psi^{r-1}\phi,v\psi^r)$, for every $r>0$. 
By uniqueness of extension, we have that $\hat\varphi$ is given by $(u,v)\mapsto (v\hat\phi,v\hat\psi)$. As above, we have that $(u,v)\hat\varphi^r=(v\hat\psi^{r-1}\hat\phi,v\hat\psi^r)$.

So $$\Fix(\hat\varphi)=\{(v\hat\phi,v)\mid v\in \Fix(\hat\psi)\}$$

 and $$\Per(\hat\varphi)=\{(v\hat\psi^{\pi_v -1}\hat\phi,v)\mid v\in \Per(\hat\psi)\},$$ where $\pi_v$ denotes the period of $v$.

\begin{proposition}
Let $\varphi$ be a type IV uniformly continuous endomorphism of $ F_n \times F_m$, with $m,n>1$. 
Then, $\text{Sing}(\hat\varphi)=\{(v\hat\phi,v)\mid v\in \text{Sing}(\hat{\psi})\}$. Consequentely, $\text{Reg}(\hat\varphi)=\{(v\hat\phi,v)\mid v\in \text{Reg}(\hat{\psi})\}$.
\end{proposition}
\begin{proof}
We start by showing that  $\text{Sing}(\hat\varphi)\subseteq \{(v\hat\phi,v)\mid v\in \text{Sing}(\hat{\psi})\}$. Take some $(v\hat\phi,v)\in (\text{Fix}(\varphi))^c$ with $v\in \text{Fix}(\hat\psi)$. Then, for every $\varepsilon>0$, the open ball of radius $\varepsilon$ centered in $(v\hat\phi,v)$ contains an element $(u_\varepsilon\phi,u_\varepsilon)\in \text{Fix}(\varphi)$, with $u_\varepsilon \in \Fix(\psi).$ Notice that $d(v,u_\varepsilon)\leq d((v\hat\phi,v),(u_\varepsilon\phi,u_\varepsilon))<\varepsilon$, thus $v\in (\Fix(\psi))^c.$

For the reverse inclusion, take some $v\in \Sing(\hat\psi)$. As above, we know that for every $\varepsilon>0$, there is some $u_\varepsilon\in B(v;\varepsilon)\cap \Fix(\psi)$. Notice that, since $\hat\phi$ is uniformly continuous, for every $\varepsilon >0$, there is some $\delta_\varepsilon$ such that, for all $x,y \in  \widehat{F_m}$ such that $d(x,y)<\delta_\varepsilon$, we have that $d(x\hat\phi,y\hat\phi)<\varepsilon$.
We want to prove that $(v\hat\phi,v)\in (\Fix(\varphi))^c$, by showing that, for every $\varepsilon >0$, the ball centered in $(v\hat\phi,v)$ with radius $\varepsilon$ contains a fixed point of $\varphi$. So, let $\varepsilon >0$ and consider
$\delta=\min\{\delta_{\varepsilon},\varepsilon\}$.
We have that $(u_\delta\hat\phi,u_\delta)\in B((v\hat\phi,v);\varepsilon)$ since, by definition of $u_\delta$, we have that $d(v,u_\delta)<\delta\leq \varepsilon$ and also, $d(v,u_\delta)<\delta_{\varepsilon}$ means that $d(v\hat\phi,u_\delta \hat\phi)<\varepsilon$.
 \end{proof}

Two infinite fixed points $(u\hat\phi,u)$ and $(v\hat\phi,v)$ are in the same $\Fix(\varphi)$-orbit if and only if there is some (finite) fixed point $(w\phi,w)$ such that $(v\hat\phi,v)=(w\phi,w)(u\hat\phi,u)$, i.e., if and only if $u$ and $v$ are in the same $\Fix(\psi)$-orbit. It is known that $\Reg(\hat\psi)$ has finitely many $\Fix(\psi)$-orbits (see \cite{[Sil13]}), so $\Reg(\hat\varphi)$ has finitely many $\Fix(\varphi)$-orbits.

Moreover, if $\psi$ is an automorphism, then every regular infinite fixed point $v$ of $\hat\psi$ must be either an attractor or a repeller (see \cite{[Sil13]}) and a singular infinite fixed point can never be an attractor nor a repeller (see \cite{[GJLL98]}).

\begin{proposition}
Let $v\in \Reg(\hat\psi)$. Then $(v\hat\phi,v)$ is an attractor if and only if $v$ is an attractor for $\hat\psi$. Moreover, if $\alpha\in \Sing(\hat\varphi)$, then $\alpha$ is not an attractor.\end{proposition}

\begin{proof} Let $(v\hat\phi,v)\in \Reg(\hat\varphi)$ such that $v$ is an attractor for $\psi$. We have that $\hat\phi$ is uniformly continuous, so, for every $\varepsilon >0$, there is some $\delta_\varepsilon>0$ such that, for every $x,y\in \hat F_m$ such that $d(x,y)<\delta_\varepsilon $ we have that $d(x\hat\phi,y\hat\phi)<\varepsilon.$ Also, let $\tau>0$ be such that $$\forall\beta\in\hat{ F_m} (d(v,\beta)<\tau\Rightarrow \lim\limits_{n\to +\infty} \beta\hat\psi^n=v).$$
We will prove that, for every $x\in \hat F_n$, $y\in B(v;\tau)$, we have $(x,y)\hat\varphi^n\to(v\hat\phi,v)$.  Let $\varepsilon>0$ and take $\delta=\min\{\varepsilon,\delta_\varepsilon\}$. Take $N\in \N$ such that for every $n>N$ we have that $d(v,y\hat\psi^n)<\delta$, which also implies that $d(v\hat\phi,y\hat\psi^n\hat\phi)<\varepsilon$. Thus, for every $n>N+1$, we have that $(x,y)\hat\varphi^n=(y\hat\psi^{n-1}\hat\phi,v\hat\psi^n)\in B((v\hat\phi,v);\varepsilon)$. Thus, $(x,y)\hat\varphi^n\to (v\hat\phi,v)$.

Hence, $(v\hat\phi,v)$ is an attractor and the set of points it attracts is given by $\hat F_n\times \mc A_{v,\hat\psi}$, where $\mc A_{v,\hat\psi}$ denotes the set of points $v$ attracts through $\psi$.

Now, we prove that, if $(w\hat\phi,w)\in \Fix(\hat\varphi)$ is an attractor, then $w$ is an attractor for $\psi$ and that completes the proof. Indeed, take $\varepsilon >0$ such that $$\forall(x,y)\in\widehat{F_n \times F_m} (d((w\hat\phi,w),(x,y))<\varepsilon\Rightarrow \lim\limits_{n\to +\infty} (x,y)\hat\varphi^n=(w\hat\phi,w)).$$
For every $y\in B(w;\varepsilon)$, we have that $d((w\hat\phi,w),(w\hat\phi,y))<\varepsilon$, so $$ \lim\limits_{n\to +\infty}(y\psi^{n-1}\phi,y\psi^n)= \lim\limits_{n\to +\infty} (w\hat\phi,y)\hat\varphi^n=(w\hat\phi,w),$$
thus $y\psi^n\to w.$
 \end{proof}

\subsection{Type VI uniformly continuous endomorphisms}
Let $\varphi:F_n\times F_m\to F_n\times F_m$ be a type VI uniformly continuous endomorphism and take $\phi\in End(F_n)$, $\psi\in End(F_m)$ such that  $ (u,v)\mapsto (u\phi,v\psi)$.

 By uniqueness of extension, we have that $\hat\varphi$ is given by $(u,v)\mapsto (u\hat\phi,v\hat\psi)$. We then have that $\Fix(\hat\varphi)=\Fix(\hat\phi)\times\Fix(\hat\psi)$,  $\Sing(\hat\varphi)=\Sing(\hat\phi)\times \Sing(\hat\psi)$ and $\Reg(\hat\varphi)=\Fix(\hat\phi)\times \Reg(\hat\psi)\cup \Reg(\hat\phi)\times\Fix(\hat\psi)$. There is no hope of finding a finiteness condition that holds in general, since if $n=2$ and $\phi$ is the identity mapping, then $\Sing(\hat\phi)$ is uncountable. In this case, both $\Reg(\hat\varphi)$ and $\Sing(\hat\varphi)$ are uncountable.

\begin{proposition}
An infinite fixed point $\alpha=(u,v)$, where $u\in \Fix(\hat\phi)$ and $v\in\Fix(\hat\psi)$ is an attractor if and only if $u$ and $v$ are attractors for $\hat\phi$ and $\hat\psi$, respectively. If, additionally, $\varphi$ is an automorphism, the same holds for repellers.
\end{proposition}
\begin{proof} Let  $\alpha=(u,v)$ be an  infinite fixed point, where $u\in \Fix(\hat\phi)$ and $v\in\Fix(\hat\psi)$. Clearly if $u\in \Fix(\hat\phi)$ and $v\in \Fix(\hat\psi)$ are attractors, then, $(u,v)\in\Fix(\hat\varphi)$ is an attractor. Indeed, in that case, there are $\varepsilon_1>0$ and $\varepsilon_2>0$ such that 
$$\forall x\in\widehat{F_n},\, \left(d(u,x)<\varepsilon_1\Rightarrow \lim\limits_{n\to +\infty} x\hat\phi^n=u\right)$$ 
and 
$$\forall y\in\widehat{ F_m},\, \left(d(v,y)<\varepsilon_2\Rightarrow \lim\limits_{n\to +\infty} y\hat\psi^n=v\right).$$
Thus, taking $\varepsilon=\min\{\varepsilon_1,\varepsilon_2\}$, we have that 
\begin{align*}
\forall (x,y)\in \widehat{ F_n\times F_m} (d((u,v),(x,y))<\varepsilon &\Rightarrow d(u,x)<\varepsilon \wedge d(v,y)<\varepsilon \\
&\Rightarrow \lim\limits_{n\to+\infty} x\hat\phi^n=u \wedge \lim\limits_{n\to+\infty} y\hat\psi^n=v\\
&\Rightarrow \lim\limits_{n\to +\infty}(x,y)\hat\varphi^n=(u,v).
\end{align*}

Conversely, suppose w.l.o.g that $u$ is not an attractor for $\hat\phi$. Then, for every $\varepsilon >0$, there is some $x_\varepsilon\in \hat{F_n}$ such that $d(u,x_\varepsilon)<\varepsilon$ but $x_\varepsilon\hat\phi^n\not\to u$. In this case, we have that, for every $\varepsilon >0$, $d((u,v),(x_\varepsilon,v))<\varepsilon$ and $(x_\varepsilon, v)\hat\varphi^n=(x_\varepsilon\hat\phi^n,v)\not\to (u,v)$.

In case, $\varphi$ is an automorphism, the repellers case is analogous.
 \end{proof}

\subsection{Type VII uniformly continuous automorphisms}
Let $\varphi:F_n\times F_m\to F_n\times F_m$ be a type VII uniformly continuous endomorphism and take $\phi:F_n\to F_m$, $\psi:F_m\to F_n$ such that  $ (u,v)\mapsto (v\psi,u\phi)$.

By uniqueness of extension, we have that $(x,y)\hat\varphi=(y\hat\psi,x\hat\phi)$. As done in the finite case, we have that an infinite fixed point must be $(x,y)\in \widehat{F_n\times F_m}$ such that $x=y\hat\psi$ and $y=x\hat\phi$ and so $x=x\hat\phi\hat\psi$ and $y=y\hat\psi\hat\phi$. So, $\Fix(\hat\varphi)\subseteq \Fix(\hat\phi\hat\psi)\times \Fix(\hat\psi\hat\phi)$. Also, for $x\in \Fix(\hat\phi\hat\psi)$, we have that $(x,x\hat\phi)\hat\varphi=(x\hat\phi\hat\psi,x\hat\phi)=(x,x\hat\phi)$. Similarly, we have that $(y\hat\psi,y)\hat\varphi=(y\hat\psi,y)$ for $y\in \Fix(\hat\psi\hat\phi)$. So, $\Fix(\hat\varphi)=\{(x,x\hat\phi)\mid x\in \Fix(\hat\phi\hat\psi)\}=\{(y\hat\psi,y)\mid y\in \Fix(\hat\psi\hat\phi)\}$. \\

Notice that, by uniqueness of extension, we have that $\hat\phi\hat\psi=\widehat{\phi\psi}$ and $\hat\psi\hat\phi=\widehat{\psi\phi}$.
\begin{proposition}
\label{regsingtip7}
Let $\varphi$ be a type VII uniformly continuous endomorphism of $ F_n \times F_m$, with $m,n>1$. 
Then, $\text{Sing}(\hat\varphi)=\{(x,x\hat\phi)\mid x\in \Sing(\hat\phi\hat\psi)\}$. Consequentely, $\text{Reg}(\hat\varphi)=\{(x,x\hat\phi)\mid x\in \Reg(\hat\phi\hat\psi)\}$.
\end{proposition}
\begin{proof}
We start by showing that  $\text{Sing}(\hat\varphi)\subseteq \{(x,x\hat\phi)\mid x\in \Sing(\hat\phi\hat\psi)\}$. Take some $(x,x\hat\phi)\in (\text{Fix}(\varphi))^c$ with $x\in \text{Fix}(\hat\phi\hat\psi)$. Then, for every $\varepsilon>0$, the open ball of radius $\varepsilon$ centered in $(x,x\hat\phi)$ contains an element $(y_\varepsilon,y_\varepsilon\phi)\in \text{Fix}(\varphi)$, with $y_\varepsilon \in \Fix(\phi\psi).$ Notice that $d(x,y_\varepsilon)\leq d((x,x\hat\phi),(y_\varepsilon,y_\varepsilon\phi))<\varepsilon$, thus $x\in (\Fix(\phi\psi))^c.$

For the reverse inclusion, take some $x\in \Sing(\hat\phi\hat\psi)$.  As above, we know that for every $\varepsilon>0$, there is some $y_\varepsilon\in B(x;\varepsilon)\cap \Fix(\phi\psi)$. Notice that, since $\hat\phi$ is uniformly continuous, for every $\varepsilon >0$, there is some $\delta_\varepsilon$ such that, for all $x,y \in  \widehat{F_n}$ such that $d(x,y)<\delta_\varepsilon$, we have that $d(x\hat\phi,y\hat\phi)<\varepsilon$.
We want to prove that $(x,x\hat\phi)\in (\Fix(\varphi))^c$, by showing that, for every $\varepsilon >0$, the ball centered in $(x,x\hat\phi)$ with radius $\varepsilon$ contains a fixed point of $\varphi$. So, let $\varepsilon >0$ and consider
$\delta=\min\{\delta_{\varepsilon},\varepsilon\}$.
We have that $(y_\delta,y_\delta\phi)\in B((x,x\hat\phi);\varepsilon)$ since, by definition of $y_\delta$, we have that $d(x,y_\delta)<\delta\leq \varepsilon$ and also, $d(x,y_\delta)<\delta_{\varepsilon}$ means that $d(x\hat\phi,y_\delta \hat\phi)<\varepsilon$.
 \end{proof}

As done in the Type IV case, observing that two infinite fixed points $(x,x\hat\phi)$ and $(y,y\hat\phi)$ are in the same $\Fix(\varphi)$-orbit if and only if $x$ and $y$ are in the same $\Fix(\phi\psi)$-orbit, we have that $\Reg(\hat\varphi)$ has finitely many $\Fix(\varphi)$-orbits.

\begin{proposition}
\label{regautos7}
Suppose that $m=n$ and both $\phi$ and $\psi$ are automorphisms. Let  $x\in \Fix(\hat\phi\hat\psi)$. Then  $x\in \Sing(\hat\phi\hat\psi)$ if and only if $x\hat\phi\in \Sing(\hat\psi\hat\phi)$. Consequentely,  $x\in \Reg(\hat\phi\hat\psi)$ if and only if $x\hat\phi\in \Reg(\hat\psi\hat\phi)$.
\end{proposition}
\begin{proof} Let $x\in \Sing(\hat\phi\hat\psi)$. Then $x\hat\phi\in \Fix(\hat\psi\hat\phi)$.  

Since $\hat\phi$ and $\hat\phi^{-1}$ are uniformly continuous, we have that for every $\varepsilon>0$, there are some $\delta_\varepsilon$, $\delta'_\varepsilon$ such that for every $z,w\in \hat F_n$, we have $$d(z,w)<\delta_\varepsilon\Rightarrow d(z\hat\phi,w\hat\phi)<\varepsilon$$
and 
$$d(z,w)<\delta'_\varepsilon\Rightarrow d(z\hat\phi^{-1},w\hat\phi^{-1})<\varepsilon.$$ Now, let $\varepsilon >0$. We want to prove that there is some $y\in \Fix(\psi\phi)\cap B(x\hat\phi;\varepsilon).$
Since $x\in \Sing(\hat\phi\hat\psi)$, then there is some $y_\varepsilon$ such that $y_\varepsilon\in \Fix(\phi\psi)\cap B(x;\delta_\varepsilon)$, so $$y_\varepsilon \phi\in \Fix(\psi\phi)\cap B(x\hat\phi;\varepsilon).$$

Now, suppose $x\in \Reg(\hat\phi\hat\psi)$. There is some $\varepsilon >0$ such that $B(x;\varepsilon)\cap \Fix(\phi\psi)=\emptyset.$ We will now prove that $B(x\hat\phi;\delta_\varepsilon')\cap \Fix(\psi\phi)=\emptyset.$ Suppose there is some $y\in B(x\hat\phi;\delta_\varepsilon')\cap \Fix(\psi\phi)$. Then, $y\phi^{-1}\in B(x;\varepsilon)$. But $$y\phi^{-1}(\phi\psi)=y\psi=y\psi\phi\phi^{-1}=y\phi^{-1},$$
so $y\phi^{-1}\in B(x;\varepsilon)\cap \Fix(\phi\psi),$ a contradiction.

 \end{proof}

\begin{proposition}
\label{atracrepel7}
Suppose that $m=n$ and both $\phi$ and $\psi$ are automorphisms. Let $x\in \Reg(\hat\phi\hat\psi)$. Then $(x,x\hat\phi)$ is an attractor (resp. repeller) if and only if $x$ is an attractor (resp. repeller)  for $\hat\phi\hat\psi$. Moreover, if $\alpha\in \Sing(\hat\varphi)$, then $\alpha$ is not an attractor nor a repeller.\end{proposition}

\begin{proof} Let $(x,x\hat\phi)\in \Reg(\hat\phi\hat\psi)$ be such that $x$ is an attractor for $\hat\phi\hat\psi$. 
Then there is some $\tau>0$ such that 
\begin{align}
\label{xatrator}
\forall\beta\in\hat{ F_n} (d(x,\beta)<\tau\Rightarrow \lim\limits_{n\to +\infty} \beta(\hat\phi\hat\psi)^n=x).
\end{align}

Also, since $\hat\phi$, $\hat\phi^{-1}$ and $\hat\psi$ are uniformly continuous, we have that for every $\varepsilon>0$, there is some $\delta_\varepsilon$, such that for every $z,w\in \hat F_n$, we have 
\begin{align*}
d(z,w)<\delta_\varepsilon \Rightarrow \quad d(z\hat\phi,w\hat\phi)<\varepsilon
\wedge  d(z\hat\phi^{-1},w\hat\phi^{-1})<\varepsilon \wedge   d(z\hat\psi,w\hat\psi)<\varepsilon.
\end{align*}

We start by proving that $x\hat\phi$ is an attractor for $\hat\psi\hat\phi$ by proving that 
\begin{align}
\label{xphiatrator}
\forall\beta\in\hat{ F_n} (d(x\hat\phi,\beta)<\delta_\tau\Rightarrow \lim\limits_{n\to +\infty} \beta(\hat\psi\hat\phi)^n=x\hat\phi).
\end{align}
Indeed, if $d(x\hat\phi,\beta)<\delta_\tau$ then $d(x\hat\phi\hat\phi^{-1},\beta\hat\phi^{-1})=d(x,\beta\hat\phi^{-1})<\tau$, so $$\lim\limits_{n\to +\infty} \beta(\hat\psi\hat\phi)^{n-1}\hat\psi =\lim\limits_{n\to +\infty} \beta\hat\phi^{-1}(\hat\phi\hat\psi)^n=x.$$ By continuity of $\hat\phi$, we have that $\lim\limits_{n\to +\infty} \beta(\hat\psi\hat\phi)^{n-1}\hat\psi\hat\phi =\lim\limits_{n\to +\infty} \beta(\hat\psi\hat\phi)^{n}=x\hat\phi.$\\

We want to prove that $(x,x\hat\phi)$ is an attractor for $\hat\varphi$. 
Let $\tau'=\min\{\tau,\delta_\tau\}$. We will prove that, for every $(z,w)\in B((x,x\hat\phi);\tau')$, we have $(z,w)\hat\varphi^n\to(x,x\hat\phi)$.  Let $\varepsilon>0$ and take $\delta=\min\{\varepsilon,\delta_\varepsilon\}$. 
By (\ref{xatrator}), we can take $N\in \N$ such that for every $n>N$ we have that 
\begin{align}
\label{ref1}
d(x,z(\hat\phi\hat\psi)^n)<\delta,\quad  \text{ and so } \quad d(x\hat\phi,z(\hat\phi\hat\psi)^n\hat\phi)<\varepsilon
\end{align} 
and by (\ref{xphiatrator}), we can take $N'\in \N$ such that for every $n>N'$ we have that 
\begin{align}
\label{ref2}
d(x\hat\phi,w(\hat\psi\hat\phi)^n)<\delta,\quad \text{ and so }  \quad d(x\hat\phi\hat\psi,w(\hat\psi\hat\phi)^n\hat\psi)=d(x,w(\hat\psi\hat\phi)^n\hat\psi)<\varepsilon.
\end{align}
 Thus, for every $n>2\max\{N,N'\}+1$, we have that, if $n$ is even, 
$(z,w)\hat\varphi^{n}=(z(\hat\phi\hat\psi)^{\frac n 2},w(\hat\psi\hat\phi)^{\frac n 2})\in B(x,x\hat\phi)$ by the first parts of both (\ref{ref1}) and (\ref{ref2}); if $n$ is odd, then $(z,w)\hat\varphi^{n}=(w(\hat\psi\hat\phi)^{\frac {n-1} 2}\hat\psi,z(\hat\phi\hat\psi)^{\frac {n-1} 2}\hat\phi)\in B(x,x\hat\phi)$ by the second parts of both (\ref{ref1}) and (\ref{ref2}).

Finally, we only have to prove that, if $(x,x\hat\phi)$ is an attractor, for $\hat\varphi$, then $x$ is an attractor for $\hat\phi\hat\psi$.  Suppose then that there is some $\tau>0$ such that 
\begin{align*}
\forall(z,w)\in\widehat{ F_n\times F_n} \left(d((z,w),(x,x\hat\phi))<\tau\Rightarrow \lim\limits_{n\to +\infty} (z,w)\hat\varphi^n=(x,x\hat\phi)\right).
\end{align*}
We will prove that 
\begin{align*}
\forall y\in\hat{ F_n} \left(d(x,y)<\tau\Rightarrow \lim\limits_{n\to +\infty} y(\hat\phi\hat\psi)^n=x\right)
\end{align*}
which is equivalent to 
\begin{align*}
\forall y\in\hat{ F_n} \left(d(x,y)<\tau\Rightarrow \forall \varepsilon >0 \;\exists N\in \N: \forall n>N \;\;\, d(y(\hat\phi\hat\psi)^n,x)<\varepsilon\right)
\end{align*}
Let $y\in \hat F_n$ be such that $d(x,y)<\tau$ and $\varepsilon>0$. We know that $d((y,x\hat\phi),(x,x\hat\phi))<\tau$, 
there is $N\in N$ such that for every $n>N$, $d((y,x\hat\phi)\hat\varphi^n,(x,x\hat\phi))<\varepsilon$. Let $n>N$. Then $$d\left((y,x\hat\phi)\hat\varphi^{2n},(x,x\hat\phi)\right)=d\left((y(\hat\phi\hat\psi)^n,x\hat\phi(\hat\psi\hat\phi)^n),(x,x\hat\phi)\right)<\varepsilon,$$
so $d(y(\hat\phi\hat\psi)^n,x)<\varepsilon.$ 

The repellers case is analogous noticing that $(x,y)\hat\varphi^{-1}=(y\hat\psi^{-1},x\hat\phi^{-1})$ and so 
\begin{align*}
\Reg(\hat\varphi^{-1})&=\left\{(x,x\hat\psi^{-1})\mid x\in \Reg(\hat\psi^{-1}\hat\phi^{-1})\right\}\\&=\left\{(x,x\hat\psi^{-1})\mid x\in \Reg\left((\hat\phi\hat\psi)^{-1}\right)\right\}.
\end{align*}
 \end{proof}

\begin{corollary}
Let $\varphi$ be a type VII automorphism of $F_n\times F_n$. Then $\alpha\in \Reg(\hat\varphi)$ is either an attractor or a repeller.
\end{corollary}
\begin{proof}  
Let $\varphi$ be a type VII automorphism of $F_n\times F_n$ and take  $\alpha\in \Reg(\hat\varphi)$. From Proposition \ref{regsingtip7}, we have that $\alpha=(x,x\hat\phi)$, for some $x\in \Reg(\hat\phi\hat\psi)$. We know, from  \cite{[Sil13]}, that $x$ is either an attractor or a repeller. The result follows directly from Proposition \ref{atracrepel7}.
 \end{proof}

\section*{Acknowledgements}
The author is grateful to Pedro Silva for fruitful discussions of these topics. The author was supported by the grant SFRH/BD/145313/2019 funded by Funda\c c\~ao para a Ci\^encia e a Tecnologia (FCT).


\begin{thebibliography}{b}
\bibitem{[Ber79]} J. Berstel, {\it Transductions and Context-free Languages}, Teubner, Stuttgart, 1979.
\bibitem{[BH92]} M. Bestvina and M. Handel, {\it Train tracks and automorphisms of free groups}, Ann. Math. 135 (1992), p. 1-51.
\bibitem{[BM16]} O. Bogopolski and O. Maslakova, {\it An algorithm for finding a basis of the fixed point subgroup of an automorphism of a free group}, International J. of Algebra and Computation 26(1) (2016), p. 29-67.
\bibitem{[Car20]} A. Carvalho, {\it On the dynamics of extensions of free-abelian times free groups endomorphisms to the completion}, arXiv:2011.05205, preprint, 2020.
\bibitem{[Car22]} A. Carvalho, {\it Eventually fixed points of endomorphisms of virtually free groups}, arXiv:2204.04543, preprint, 2022.
\bibitem{[CS09b]} J. Cassaigne and P. V. Silva, {\it Infinite words and confluent rewriting systems: endomorphism extensions}, Internat. J. Algebra Comput., 19.4:443–490, 2009.
\bibitem{[CH10]} L. Ciobanu and A. Houcine, {\it The monomorphism problem in free groups}, Archiv der Mathematik, 94(5):423–434, 2010.
\bibitem{[Coo87]} D. Cooper, {\it Automorphisms of free groups have finitely generated fixed point sets}, J. Algebra 111 (1987), p. 453-456.
\bibitem{[Day09]} M.B. Day, {\it Peak reduction and finite presentations for automorphism groups of right-angled artin groups}, Geom. Topol., 13:817–855, 2009.
\bibitem{[Day14]} M.B. Day, {\it Full-featured peak reduction in right-angled artin groups}, Algebr. Geom. Topol., 14(3):1677–1743, 2014.
\bibitem{[DV13]} J. Delgado and E. Ventura, {\it Algorithmic problems for free-abelian times free groups}, J. Algebra 263(1) (2013), p. 256-283.
\bibitem{[DJK16]} V. Diekert, A. Je\.z and M. Kufleitner, {\it Solutions of Word Equations Over Partially Commutative Structures}, 43rd International Colloquium on Automata, Languages, and Programming (ICALP), 55 (2016) 127:1-127:14.
\bibitem{[GJLL98]} D. Gaboriau, A. Jaeger, G. Levitt and M. Lustig, {\it An index for counting fixed points of automorphisms of free groups}, Duke Math. J. 93 (1998), p. 425-452. 
\bibitem{[Ger87]} S. M. Gersten, {\it Fixed points of automorphisms of free groups}, Adv. Math. 64
(1987), p. 51-85.  
\bibitem{[How54]} A.G.Howson, {\it On the intersection of finitely generated free groups}, London Math. Soc., 29:428–434, 1954.
\bibitem{[LL08]} G. Levitt and M. Lustig, {\it Automorphisms of free groups have asymptotically periodic dynamics}, J. Reine Angew. Math. 619 (2008), p. 1-36.
\bibitem{[Lot83]}  M. Lothaire, {\it Combinatorics on Words}, Addison-Wesley, Reading, 1983.
\bibitem{[Mak82]} Makanin, {\it G. equations in free groups (russian)}, arxiv:1806.09560, Izv. Akad. Nauk SSSR Ser. Mat,
46:1190–1273, 1982.
\bibitem{[MS18]} F. Matucci and P. V. Silva, {\it Extensions of automorphisms of self-similar groups}, arxiv:1806.09560, preprint (CMUP 2018-8).
\bibitem{[Mih58]}  K. A. Mihailova, {\it The occurrence problem for direct products of groups}, Dokl. Acad.
Nauk SSRR, 119:1103–1105, 1958.
\bibitem{[Mut21]}  J. P. Mutanguha, {\it Constructing stable images}, preprint, available at https://mutanguha.com/pdfs/relimmalgo.pdf.
\bibitem{[MS02]} A. G. Myasnikov and  V. Shpilrain, {\it Automorphic orbits in free groups},   J. Algebra 269 (2002), p. 18-27.
\bibitem{[Pau89]} F. Paulin, {\it Points fixes d’automorphismes de groupes hyperboliques}, Ann. Inst.
Fourier 39 (1989), p. 651-662.
\bibitem{[RSS13]} E. Rodaro, P. V. Silva and M. Sykiotis, {\it Fixed points of endomorphisms of graph
groups}, J. Group Theory 16(4) (2013), p. 573-583.  
\bibitem{[Sil13]} P. V. Silva, {\it Fixed points of endomorphisms of virtually free groups}, Pacific J.
Math. 263(1) (2013), p. 207-240.
\bibitem{[Whi36]} J. H. C. Whitehead, {\it On equivalent sets of elements in a free group}, Ann. of Math., 37(4):782–800, 1936.
\bibitem{[SW70]} Stephen Willard, {\it General Topology}, Addison-Wesley, Reading, Mass., 1970.  
%
\end{thebibliography}
\end{document}